\documentclass[final,leqno,onefignum,onetabnum]{siamltex1213}

\usepackage{soul} 

\usepackage{graphicx,color,amsbsy}
\usepackage{amsfonts}
\usepackage{amsmath}
\usepackage{latexsym}
\usepackage{amssymb}
\usepackage{amscd}
\usepackage{enumerate}
\usepackage{comment}
\usepackage{multicol}
\usepackage{multirow}
\usepackage{color}
\usepackage{relsize}
\usepackage{xspace}
\usepackage{tikz}
\usepackage{pgflibraryarrows}

\title{Comparison of adaptive multiresolution and adaptive mesh refinement applied
to simulations of the compressible Euler equations\footnote{An early version of this paper has appeared previously in {\em ESAIM Proceedings}, 16:181--194, 2009 [21].}}

\author{Ralf Deiterding\thanks{Aerodynamics and Flight Mechanics Research Group, University of Southampton, Highfield Campus, Southampton SO17 1BJ, United Kingdom, \email{r.deiterding@soton.ac.uk}} 
\and Margarete O. Domingues \thanks{Laborat\'orio Associado de Computa\c c\~ao e Matem\'atica Aplicada (LAC), Coordenadoria dos 
Laborat\'orios Associados (CTE), Instituto Nacional de Pesquisas Espaciais (INPE),
Av. dos Astronautas 1758, 12227-010 S\~ao Jos\'e dos Campos, S\~ao Paulo, Brazil, \email{margarete.domingues@inpe.br}}
\and S\^onia M. Gomes\thanks{Universidade Estadual de Campinas (Unicamp), IMECC, Rua S\'ergio Buarque de Holanda, 651, Cidade Universit\'aria Zeferino Vaz, 13083-859, Campinas – SP – Brazil,  \email{soniag@ime.unicamp.br} }
\and~Kai~Schneider\thanks{M2P2-CNRS \& Centre de Math\'ematiques et d'Informatique (CMI), Universit\'e d'Aix-Marseille, 
39 rue F. Joliot-Curie, 13453 Marseille Cedex 13, France,  \email{kschneid@cmi.univ-mrs.fr}}
}

\newcommand{\correction}[1]{{#1}}
\newcommand{\Rplus}{\protect\hspace{-.1em}\protect\raisebox{.35ex}{\smaller{\smaller\textbf{+}}}}
\newcommand{\Cpp}{\mbox{C\Rplus\Rplus}\xspace}
\begin{document}


\maketitle
\slugger{mms}{xxxx}{xx}{x}{x--x}

\graphicspath{{FIG/}}

\begin{abstract} 
We present a detailed comparison between two adaptive numerical approaches to solve partial differential equations (PDEs), adaptive multiresolution (MR) and adaptive mesh refinement (AMR).
Both discretizations are based on finite volumes in space with second order shock-capturing,
and explicit time integration either with or without local time-stepping.
The two methods are benchmarked for the compressible Euler equations in Cartesian geometry.
As test cases a 2D Riemann problem, Lax-Liu $\#6$, and a 3D ellipsoidally expanding shock wave have been chosen.
We compare and assess their computational efficiency in terms of CPU time and memory requirements.
We evaluate the accuracy by comparing the results of the adaptive computations with those obtained with the corresponding FV scheme using a regular fine mesh.
We find that both approaches yield similar trends for CPU time compression for increasing number of refinement levels.
MR exhibits more efficient memory compression than AMR and shows slightly enhanced convergence; however, a larger absolute overhead is measured for the tested codes.
\end{abstract}

\begin{keywords}
adaptive numerical methods, conservation laws, Euler equations, multiresolution, mesh refinement, local time stepping
\end{keywords}

\begin{AMS}65M50, 65Y20, 76M12\end{AMS}

\pagestyle{myheadings}
\thispagestyle{plain}
\markboth{R. DEITERDING, M. O. DOMINGUES, S. M. GOMES, AND K. SCHNEIDER }{COMPARISON OF ADAPTIVE MULTIRESOLUTION AND AMR}

\section{Introduction}

Adaptive discretization methods for solving nonlinear PDEs have a long tradition in scientific computing, see \textit{e.g.}, \cite{Brand-77}. 
They are motivated by the computational complexity of real world 
problems which often involve a multitude of active spatial and temporal scales. 
Adaptivity can be understood in the sense that the computational 
effort is concentrated at locations and time instants where it is necessary to ensure a given 
numerical accuracy, while efforts may be significantly
reduced elsewhere.
Typical applications are combustion problems with thin chemical reaction zones, fluid and plasma turbulence showing self-organization into coherent vortices, and 
more generally, most kinds of interface and boundary layer type problems.

One of the essential ingredients of fully adaptive schemes is a reliable error 
estimator for the solution. 
It can be provided by Richardson extrapolation, auxiliary solutions of adjoint 
problems \cite{BR01}, gradient based approaches or wavelet coefficients
\cite{RSTB03,CohenKaberMullerPostel:2003,Cohen:2000,Muller:2003,DominguesGomesRousselSchneiderESAIM:2011,Kolomenskiy:2015,Rossinellietal:2015}.
However, there is a price for adaptivity, since the computational cost per unknown may increase 
significantly.
Hence, an adaptive method can only be efficient when the data compression is sufficiently 
large  to compensate the additional computational cost, and this is problem dependent. 
To perform efficient adaptive simulations there is an effort that has to be made to  
design algorithms and data structures, the latter are usually based on graded trees, hash-tables, linked lists or multi-domains.
%
%

In the present paper, adaptive computations of 2D and 3D  compressible 
Euler equations are presented. 
The goal is to compare, for the same prescribed accuracy, the efficiency in terms of 
CPU time and memory compression of two  approaches: the adaptive multiresolution (MR) method and the adaptive mesh 
refinement (AMR) method.
We consider global time stepping and scale-dependent local time stepping. 
Both methods are based on explicit finite volume (FV) discretizations on  adaptive meshes, with  schemes in space and time of second order accuracy. 
Our interest in comparing MR and AMR methods in the current paper can be seen as a  step towards detailed  benchmarking, 
which has been started with a preliminary 2D case study in \cite{DeiterdingDominguesGomesRousselSchneiderESIAM:2009}.
Beside considering 3D results, we use an improved version of the MR  Carmen code where the memory allocation has been changed and the computational efficiency of the underlying FV scheme has been increased.
A detailed analysis of CPU time and memory compression including error assessment allows for a sound evaluation
of the MR and AMR methods.
%
%

The block-structured AMR technique for hyperbolic PDEs has been pioneered by Berger and Oliger \cite{Berger-Oliger-84}. 
While the first approach utilized rotated refinement meshes, AMR denotes today 
especially the variant of Berger and Collela \cite{Berger-Collela-88} that 
allows only refinement patches aligned to the coarse mesh. The 
efficiency of this algorithm, in particular for instationary supersonic gas dynamical problems,
has been demonstrated in Bell {\em et al.} \cite{Bell-Berger-Saltzman-94}. Several
implementations of the AMR method for single processor computers 
\cite{Berger-LeVeque-98,Friedel-Grauer-Marliani-97} and 
parallel systems \cite{Bell-Berger-Saltzman-94,Kohn-Baden-95,Rendleman-etal-00,AMROC} have been presented; variants
utilizing simplified data structures have also been proposed \cite{MacNeice-etal-00}. For an overview of the AMR method and its implementation we refer to \cite{DeiterdingESAIM:2011}.

Multiresolution techniques have been developed after AMR,  and became popular for 
hyperbolic conservation laws going back to the 
seminal work of Harten \cite{Harten:1995} in the context of FV schemes
and cell-average MR analysis.
Starting point is an FV scheme for hyperbolic conservation laws on a regular mesh. 
Subsequently, a discrete multiresolution analysis is used to avoid expensive flux computations in smooth regions, first without reducing memory requirements,
\textit{e.g.}, for 1D  hyperbolic conservation laws \cite{Harten:1995,Bihari:1996}, 2D hyperbolic conservation laws \cite{BihariHarten:1997}, 2D compressible Euler equations \cite{ChiavassaDonat:2001}, 2D hyperbolic conservation laws with curvilinear patches \cite{DGM:2001} and unstructured meshes 
\cite{AbgrallHarten:1998,CDKP:2000}. 
A fully adaptive version, still in the context of 1D and 2D hyperbolic conservation laws, has been developed to reduce also memory requirements 
\cite{Gottschlich:1999,Kaibara:2000,CohenKaberMullerPostel:2003}. 
This algorithm has been extended to 3D and to parabolic PDEs \cite{RSTB03,RousselSchneider:2005}, and more recently to self-adaptive global and local time-stepping by M\"uller and Stiriba and ourselves in
\cite{MullerStiriba:2007,DGRS:2008,DGRS:2009,DRS:2009}. 
Therewith the solution is represented and computed on a dynamically evolving automatically adapted mesh. Different strategies have been proposed to evaluate the flux without requiring full knowledge of fine mesh cell-average values.
The MR approach has also been used in other discretization contexts. 
For instance, the Sparse Point Representation (SPR) method is the first fully adaptive MR scheme, introduced in \cite{Holmstrom:1997,Holmstrom:1999} in the context
of finite differences and point-value MR analysis, leading to both CPU time and memory reduction. 
This approach has  been  explored in applications of the SPR method in \cite{DominguesGomesDiaz:2003}, and more recently  in \cite{RossinelliHejazialhosseiniSpampinatoKoumoutsakos:2011} for parallel  implementation of the MR method for FV discretizations. 
%
%
For comprehensive literature about the subject, we refer to the books \cite{Cohen:2000,Muller:2003} and our review papers \cite{SchneiderVasilyev:2010,DominguesGomesRousselSchneiderESAIM:2011}.



%
The outline of the paper is the following:
first, we sketch the set of compressible Euler equations together with their finite 
volume discretization in space and an explicit scheme in time.   
Then, we briefly summarize the MR and AMR strategies.
In the main part, we show and discuss the numerical results for two test problems, one in 2D and one in 3D. Final remarks and conclusions based on the performed computations are drawn in Section~\ref{conclusions}.

\section{Numerical methods}\label{sec:nummethods}

For the study of the present paper we consider FV discretizations of the 3D compressible Euler equations given in the conservation form
\begin{equation}
\partial_t\bf q + \nabla \cdot {\bf f} ({\bf q})= 0, 
\label{conservation}
\end{equation}
with ${\bf q} = \left(\rho,  \rho {\bf v},  \rho e \right)^T$, where $\rho=\rho({\bf x},t)$ is the density, 
${\bf v}={\bf v}({\bf x},t)$ is the  velocity vector with components $(v_1,v_2,v_3)$, and
$e=e({\bf x},t)$ is the energy per unit of mass. All variables  are functions of time $t$ and   position ${\bf x}=(x_1,x_2,x_3)$. 
The flux function ${\bf f}=(f_{1},f_{2},f_{3})^{T}$ is given by
\begin{align}
f_{1} & =\left(\begin{array}{c}
\rho v_{1}\\
\rho v_{1}^{2}+p\\
\rho v_{1}v_{2}\\
\rho v_{1}v_{3}\\
(\rho e+p)v_{1}\\
\end{array}\right),\;\; 
f_{2}=\left(\begin{array}{c}
\rho v_{2}\\
\rho v_{1}v_{2}\\
\rho v_{2}^{2}+p\\
\rho v_{2}v_{3}\\
(\rho e+p)v_{2}\\
\end{array}\right),\;\;
f_{3}=\left(\begin{array}{c}
\rho v_{3}\\
\rho v_{1}v_{3}\\
\rho v_{2}v_{3}\\
\rho v_{3}^{2}+p\\
(\rho e+p)v_{3}\\
\end{array}\right),
\end{align}

\noindent where the pressure $p=p({\bf x},t)$ satisfies the  ideal gas constitutive relation
\begin{equation}
p = \rho R T = \left( \gamma - 1 \right) \rho \left( e - \frac{|{\bf v}|^2}{2}\right),
\end{equation}
where $T=T({\bf x},t)$ is the temperature, $\gamma=1.4$  the specific heat ratio and $R$ a specific gas constant.
In the 2D formulation, the $3^{rd}$ components of the vectors vanish and the variables only depend on $x_1$ and $x_2$.


As reference discretization for equations in the conservation form (\ref{conservation}), we consider the numerical solution represented
by the vector ${\bf Q}(t)$ of the approximated cell averages \vspace*{-0.2cm}
\begin{equation}
{\bf Q}_{i,j,k}(t)=\frac{1}{|\Omega_{i,j,k}|}\int_{\Omega_{i,j,k}}{\bf q}(\bf x,t) \; d{\bf x}
\end{equation}
on cells $\Omega_{i,j,k}$ of a uniform mesh of the computation domain \(\Omega\). For  space discretization, a FV  method is 
chosen, which results in a system of ordinary differential equations of the form \vspace*{-0.2cm}
\begin{equation}
\frac{d{\bf Q}}{dt}=-{\bf F}({\bf Q}),\label{eq:ODE}
\end{equation}
where ${\bf F}({\bf Q})$ denotes the vector of numerical flux function differences with respect to each cell. 
In all numerical schemes throughout this paper enhanced AUSM-type numerical flux functions with comparable second order accurate reconstruction and flux limiting are used. 
For all MR simulations, including the corresponding reference FV solution, a second order MUSCL scheme with an AUSM+ flux vector splitting scheme \cite{Liou:1996} and  van Albada limiter is applied. 
In all AMR computations, also including the corresponding reference FV solution, 
a standard unsplit shock-capturing MUSCL scheme with Minmod limiter and AUSMDV flux-vector splitting is adopted \cite{Wada-Liou-97}.  

For  time integration, approximate solutions  ${\bf Q}{}^{n}$ at a sequence of time instants $t^{n}$ are obtained 
using explicit ordinary differential equation solvers. 
Here, an explicit second order Runge-Kutta (RK2) scheme is used for MR,  and the MUSCL-Hancock approach \cite{Deiterding-PhDThesis,Toro} is 
used for  AMR. While both approaches correspond to using a second order accurate explicit midpoint rule for temporal integration, the key difference
of the MUSCL-Hancock method is that it utilizes the {\em exact} flux function in the intermediate time step instead of the AUSM-type numerical flux. This makes this method
slightly less accurate but computationally cheaper.

In summary,  six numerical discretizations are considered. 
For both MR and AMR, FV reference solutions are computed using the corresponding numerical scheme. 
In addition, we perform in both cases either global or local time stepping.
In the following subsections we briefly describe first the MR and then the AMR method.

\subsection{Adaptive multiresolution method}

The adaptive MR scheme belongs to a class of adaptive hybrid methods which are formed by two basic parts: the operational part and the representation part.
The operational part consists of an accurate and stable discretization of the partial differential operators.
In the representation part,  wavelet tools are employed for the MR analysis of the discrete information.

\begin{figure}[t]
\centering
\begin{tabular}{lr}
\includegraphics[height=3.5cm]{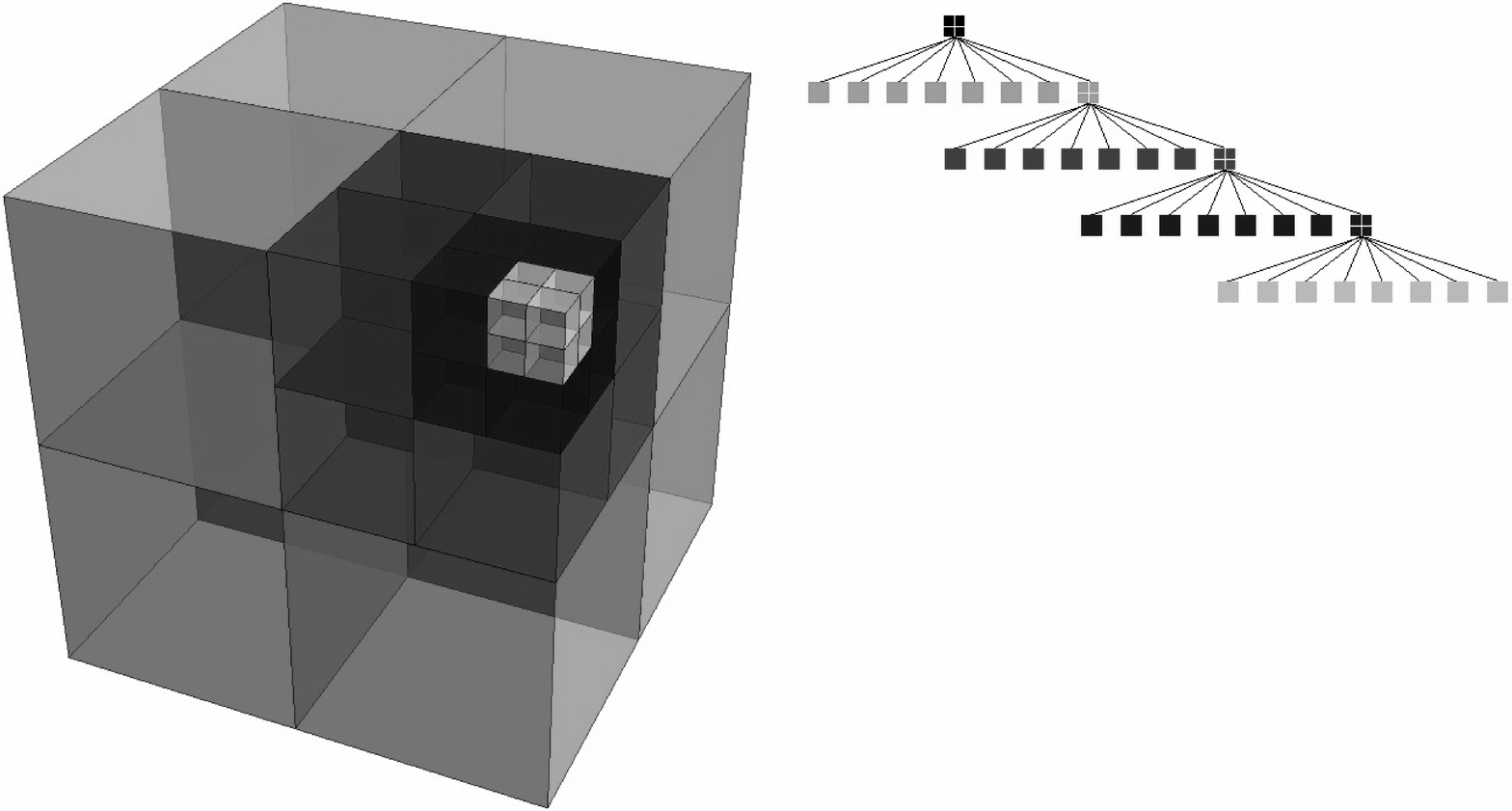} &
\includegraphics[width=3cm]{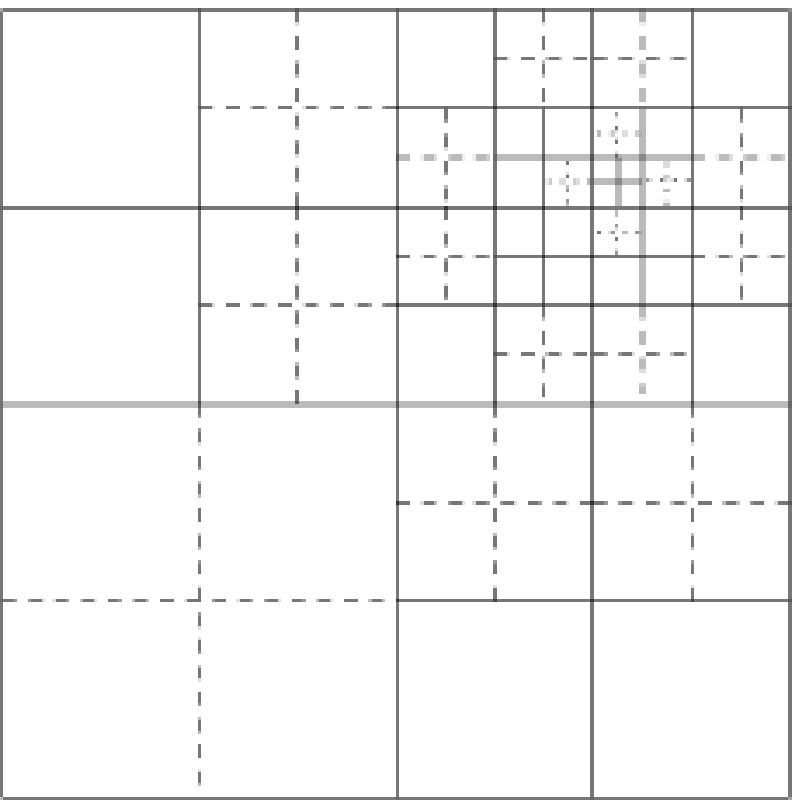} \\
& \\
\end{tabular}
\caption{Leaves of an adaptive 3D MR mesh, left. 
Sketch of the corresponding octree structure, middle. Leaves ({\it plain}) and virtual leaves ({\it dashed}) of an adaptive 2D MR mesh, right.}
\label{fig:tree}
\end{figure}

The principle in MR methods is the transformation of the cell averages given by the FV discretization into a multiscale representation.  
A hierarchy of nested meshes endowed with projection and prediction operators to perform  the inter-level transformations are the main building blocks \cite{Harten:1995,Harten:1996}. 
The numerical solution at the highest resolution level is transformed  into a  set of coarser scale approximations  plus a series of prediction errors corresponding to wavelet coefficients. These coefficients describe the difference between subsequent resolutions. 
The main idea is then to use the decay of the wavelet coefficients to estimate the local regularity of the solution \cite{Cohen:2000,Muller:2003}.
In regions where the solution is smooth these coefficients are small, while they have large magnitude in regions of steep gradients or discontinuities.
An adaptive MR representation of the numerical solution can then be obtained by a compact multiscale representation which is constructed  by removing the wavelet coefficients whose magnitude is smaller than a chosen threshold \cite{RSTB03,CohenKaberMullerPostel:2003}.  
A natural way to store the compact MR representation is to use a tree data structure. Locally refined adaptive meshes create incomplete trees, cf. Fig.~\ref{fig:tree}.
The adaptive computations are performed on the leaves of the tree, defined as nodes without children, where the cell averages are stored.
In this procedure, gradedness of the tree is imposed, \textit{i.e.}, only one level of difference is permitted between subsequent neighbors.
%
For the complete time evolution cycle, three basic operations are considered:

{\bf Refinement:} 
The solution may change in time,
the adaptive MR mesh at $t^{n}$ may not be sufficient at the next time
step. 
Thus, a preventive action is necessary 
to account for possible translation or creation of finer structures in
the solution between subsequent time steps.
Before performing the time evolution, the solution is interpolated onto an extended mesh.  For this the adaptive mesh at $t^{n}$ is refined by one level while maintaining the gradedness.

{\bf Time evolution:}
First, to ensure conservation of the flux computations between different levels it is necessary to add virtual leaves, as illustrated in Figure~\ref{fig:tree}, right. %
The time evolution operation is applied only on the leaves of the extended mesh, but not on the virtual ones \cite{RSTB03}.
There is a possibility to save further CPU time by evolving the solution with level dependent time steps, instead of a global time step.
This MRLT scheme with a given time step $\Delta t^n$ at the finest level $L$,
evolves the cells at coarser levels $0\leq\ell<L$ with larger time steps
$\Delta t_\ell^n =2^{L-\ell}\Delta t^n$. 
Required missing values are interpolated in time at intermediate time steps.

{\bf Coarsening:} 
The regions of smoothness of the solution can change after the time evolution. Hence, the adaptive MR mesh must be checked if a coarser mesh is sufficient for an accurate representation of the computed solution.
A multiresolution analysis is performed and the thresholding of the wavelet coefficients determines the cells where the mesh can be coarsened.

\smallskip

The cell-average MR analysis used in the present paper
corresponds to a third order prediction operator based on  quadratic polynomial interpolation of the cell averages. 
The adaptive MR scheme and its MRLT version, cf. \cite{DGRS:2008,DominguesGomesRousselSchneiderESAIM:2011}, are implemented in the code Carmen\footnote{\correction{The Carmen code is open access and available at \\ \url{https://github.com/waveletApplications/carmen.git}.}} originally developed by  Roussel \cite{RSTB03}.
Here an optimized version of the Carmen code is used where the runtime has been improved for both the MR and the FV method. 
The implementation is in \correction{\Cpp} throughout and consists of approximately $20,000$ lines of code (LOC) in total.

\subsection{Adaptive mesh refinement method}\label{sec:amr}

\begin{figure}[t]
\begin{center}
\parbox{7cm}{\centering\footnotesize{\tt
\begin{tabbing}
\hspace*{0.0cm} \= \hspace*{0.3cm} \= \hspace*{0.3cm} \= \hspace*{0.3cm} \= \hspace*{0.3cm} \kill
AdvanceLevel($l$) \\[0.2ex]
\> Repeat $r_l$ times \\
\> \> Set ghost cells of ${\bf Q}^l(t)$ \\
\> \> If time to remesh \\
\> \> \> Remesh($l$)\\
\> \> UpdateLevel($l$) \\
\> \> If level $l+1$ exists \\
\> \> \> Set ghost cells of ${\bf Q}^l(t+\Delta t_l)$ \\
\> \> \> AdvanceLevel($l+1$) \\
\> \> \> Average ${\bf Q}^{l+1}(t+\Delta t_l)$ onto \\
\> \> \> \>  ${\bf Q}^l(t+\Delta t_l)$\\
\> \> \> Flux correction of ${\bf Q}^l(t+\Delta t_l)$ \\
\> \> $t:=t+\Delta t_l$
\end{tabbing}}}\hspace*{0.5cm}
\parbox{7cm}{\centering\footnotesize
\includegraphics[width=6.5cm]{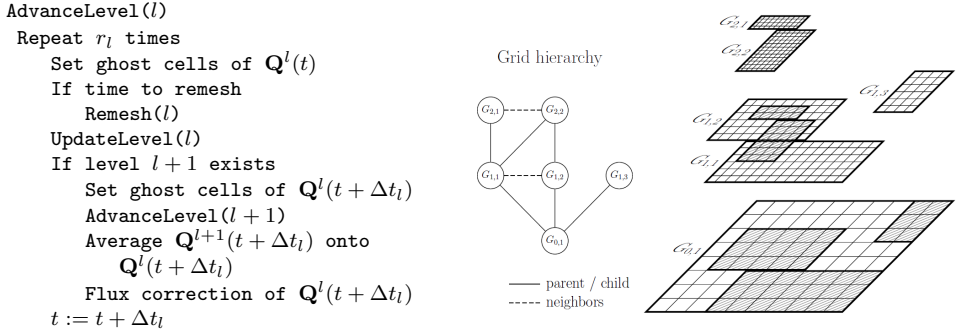}}
\end{center}
\caption{\label{pic:hier} Recursive AMR algorithm in 2D, and typical hierarchy of rectangular submeshes.}
\end{figure}

The AMR method \cite{Berger-Oliger-84,Berger-Collela-88,Bell-Berger-Saltzman-94} follows a patch-oriented refinement approach, where
non-overlapping rectangular submeshes 
$G_{l,m}$ define the domain  $G_l := \bigcup_{m=1}^{M_l} G_{l,m}$ of an entire level $l=0,\dots, L$.
As the construction of refinement proceeds recursively, a hierarchy of submeshes 
successively contained within the next coarser level domain is created, as illustrated in Figure~\ref{pic:hier}. 
Note that
the recursive nature of the algorithm allows only the addition of one new level in each refinement operation.  The patch-based approach does not require  special
coarsening operations; submeshes are simply removed from the hierarchy. The coarsest possible resolution is thereby
restricted to the level 0 mesh. Usually, it is assumed that all mesh widths on level $l$ are $r_l$-times finer 
than on the level $l-1$, \textit{i.e.},  $\Delta t_l=\Delta t_{l-1}/r_l$ and $\Delta x_{n,l}=\Delta x_{n,l-1}/r_l$, with 
$r_l\in\mathbb{N}, r_l\ge 2$ for $l>0$ and $r_0=1$, which ensures that a time-explicit FV scheme remains stable under a CFL-type condition on all levels of the hierarchy. 

The numerical update is applied on the level $l$ by calling a single-mesh routine implementing the 
FV scheme in a loop over all the submeshes $G_{l,m}$. The regularity of the input 
data allows a straightforward implementation of the scheme and further permits optimization to take 
advantage of high-level caches, pipelining, etc. New refinement meshes are initialized 
by interpolating the vector of conservative quantities ${\bf Q}$ from the next coarser level. However, data in cells already 
refined is copied directly from the previous refinement patches. {\it Ghost} (also know as {\it halo}) cells around each patch 
are used to decouple the submeshes computationally. Ghost cells outside of the root domain $G_0$ are 
used to implement physical boundary conditions. Ghost cells in $G_l$ have a unique interior cell 
analogue and are set by copying the data value from the patch where the interior cell is contained 
(synchronization). For $l>0$, internal boundaries can also be used. If recursive time step refinement
is employed, ghost cells at the internal refinement boundaries on the level $l$ are set by time-space interpolation 
from the two previously calculated time steps of level $l-1$. Otherwise, spatial interpolation from the
level $l-1$ is sufficient.

The characteristic of the AMR algorithm is that refinement patches overlay coarser mesh data structures,
instead of being embedded, again avoiding data fragmentation. Values of cells covered by finer 
submeshes are subsequently overwritten by averaged fine mesh values, which, in general, would
lead to a loss of conservation on the coarser mesh. A remedy to this problem is to replace 
the coarse mesh numerical fluxes at refinement boundaries with the sum of fine mesh fluxes along 
the corresponding coarse cell boundary. Details about this flux correction can be found in 
\cite{Berger-Collela-88,Deiterding-PhDThesis,Deiterding-03}.
The basic recursive AMR algorithm is formulated in Fig.~\ref{pic:hier} (left). New refinement meshes on all the higher levels are
created by calling {\tt Remesh()} at a given level $l$. The level $l$ by itself is not modified. 
To consider the nesting of the level domains already in the mesh generation, {\tt Remesh()} 
starts at the highest level to be refined and proceeds down to $l+1$. After evaluating the
refinement indicators and flagging cells for refinement, a special clustering algorithm 
\cite{Bell-Berger-Saltzman-94} is used to create new refinement patches until the ratio between flagged and all cells 
in every new submesh is above a prescribed threshold $0<\eta_{tol}\le 1$. 

In the present paper, all the AMR computations have been carried out using the AMROC (Adaptive Mesh Refinement in 
Object-oriented \correction{\Cpp}) framework \cite{AMROC,VTF}\footnote{\correction{The latest open access release of AMROC is available at \url{http://www.vtf.website}.}}. 
At the present time, the AMR core of AMROC consists of 
approximately $46,000$ LOC in \correction{\Cpp} and approximately $6,000$ LOC for visualization 
and data conversion. Similarly to the AMR inter-level transfer operations (interpolation, averaging), the 
employed FV update routine is coded in Fortran-77 and all the codes are compiled with
standard compiler optimizations (-O3 with loop unrolling, inlining, etc.) using the GNU compiler suite on the benchmark system with Intel-i7-$2.2\,\mathrm{GHz}$ quad-core processor.  
Although AMROC permits large-scale MPI-parallel AMR computations \cite{Deiterding-08c,Pantano-Deiterding-05}, the present investigation 
uses only the {\it serial} algorithm of the software. 

The original recursively adaptive AMR method with time step refinement is denoted here as AMRLT. In order to provide a comparison 
to the MR method, we also employ the same implementation \correction{without local time stepping} under the name AMR. Note that the performance of this variant could 
be improved by skipping some of the adaptation steps provided that refinement coverage is correspondingly enlarged. Yet, such 
optimizations have not been investigated here. 
As refinement indicators, scaled gradient criteria of the form 
\begin{equation} 
|w(Q_{i+\alpha_i,j+\alpha_j,k+\alpha_k})-w(Q_{i,j,k})|  > \epsilon_w, 
\end{equation}
with $\alpha_i, \alpha_j, \alpha_k \in \{0,1\}$, are applied to density and pressure ($w=\rho$ and $w=p$). As it is common practice \cite{Berger-Collela-88}, a
layer of one additional cell width is also tagged for the refinement around each refinement flag to ensure that the flagged feature 
does not leave the refinement region during the next time step. 
Furthermore, AMROC allows for the additional application of a heuristic local error indicator based on a Richardson estimation 
\cite{Berger-Collela-88,Berger-Oliger-84}. For the test problems of Section~\ref{numerical}, however, 
scaled gradient and Richardson error estimation criteria were found to give virtually identical mesh refinement and the 
benchmarked computations only utilized the former.

\section{Error assessment}

The goal of the adaptive MR and AMR computations is to obtain the solution with a significant gain in CPU time and  memory while preserving the accuracy of the corresponding FV scheme on the regular finest mesh. 
%
To quantify the accuracy of the adaptive simulations the discrete $L_1$-error is computed.
For that a reference  FV solution on a finer mesh  restricted onto the current uniform mesh level $L$ is used.

For the MR  methods, with or without local time stepping ($a=\mathrm{MR}$ or MRLT), the adaptive solution is recursively projected up to the desired finest uniform mesh of level $L$ with a step size $\Delta x_{L}$. 
The goal is to obtain $\tilde{\bf Q}_{(i,j)}^a$ using third order cell-average interpolation.   
The discrete error is then evaluated on the domain $\Omega$ as 
\begin{equation}
L_1({\bf Q}) = \sum_{i,j,k} |\tilde{\bf Q}_{i,j,k}^a-{\bf Q}_{i,j,k}^r| \Delta x_{L}^2, 
\end{equation}
where ${\bf Q}^r_{i,j,k}$ denotes the projection of the reference solution down to the desired level $L$.
For the two  AMR cases, considering again either global or local time stepping, the error is evaluated as the sum of the $L_1$-error computed on the domain $\Omega_l$ without higher refinement, \textit{i.e.},   \vspace*{-0.2cm}
\begin{equation}
L_1({\bf Q}) = L_1(\Delta x_L,\Omega_L) + \sum_{l=0}^{L-1} L_1(\Delta x_l,\Omega_l\setminus \Omega_{l+1}), 
\end{equation}
where \vspace*{-0.2cm}
\begin{equation}
L_1(\Delta x,\Omega) = \sum_{i,j,k} |{\bf Q}_{i,j,k}^a-{\bf Q}_{i,j,k}^r| \Delta x^2, 
\end{equation}
denotes the $L_1$-norm on the domain $\Omega$.
The projection of the reference solution from the finer uniform mesh down to the desired mesh with step size $\Delta x$ is denoted by ${\bf Q}^r_{i,j,k}$. 
The index $a$ stands for either AMR or AMRLT.

To evaluate the performance of the adaptive codes, \textit{i.e.},  CPU, memory and mesh compression rates, including accuracy perturbation, and overhead,  different measures are introduced.
The  CPU time compression rate is defined  as the ratio between the CPU time required to compute the solution at the final instant using the adaptive method $a$ and the one required to compute the same solution using the fine-mesh FV method  \vspace*{-0.1cm}
\begin{equation}
\mbox{\it CPU time compression rate}= \frac{CPU^a}{CPU^{FV}},
\end{equation}
The ratio of the average memory requirement and the number of cells $N_C$ of the finest uniform mesh defines the memory compression rate in the adaptive computations,
\begin{equation}
\mbox{\it memory compression rate} = \frac{\frac{1}{N_I^a} \sum_{n=1}^{N_I^a} {\cal C}^{a,n}}{N_C},
\end{equation}
where  $N_{I}^a$ is the number of performed time steps, and ${\cal C}^{a,n}$ denotes the sum of cells of the entire hierarchy at $t=t^n$. 
Similarly,  {\it mesh compression} can be defined as the ratio of the average leaf requirement (\textit{i.e.}, the average number of required degrees of freedom to represent the numerical solution) and the  number of cells of the finest uniform mesh  \vspace*{-0.4cm}
\begin{equation}
\mbox{\it mesh compression rate}=\frac{\frac{1}{N_I^a} \sum_{n=1}^{N_I^a} {\cal L}^{a,n}}{N_C}.
\end{equation}
In the adaptive codes, typically  a perturbation error  is introduced.
This results in larger errors of the adaptive solution compared to the FV solution on the finest uniform mesh. 
Suitable thresholding guarantees that this perturbation does indeed not deteriorate the convergence order of the underlying  FV scheme. 
The accuracy perturbation is measured by the relative difference between the error given by the reference FV scheme and the error introduced by the adaptive method,
\begin{equation}
\mbox{\it accuracy perturbation}=\frac{| L_1^{FV}({\bf Q})- L_1^{a}({\bf Q})|}{L_1^{FV}({\bf Q})}.
\end{equation} 
Evolving the solution in each cell of the computational domain using adaptive computations is expected to be more expensive in terms of CPU time than using  the  reference scheme on the uniform mesh.  
Hence, an essential  question is to know whether the reduction in the number of degrees of freedom required to represent the adaptive numerical solution compensates the additional computational overhead induced by the adaptive algorithm. 
The overhead of the adaptive computations can be evaluated by 
\begin{equation} 
\Gamma^{a}= \frac{CPU^a} {\sum_{n=1}^{N_I^a} {\cal L}^{a,n}},
\end{equation}
which denotes the average CPU time per leaf spent to evolve the solution in the computational domain.  
Therefore, the ratio $\Gamma^{a}/\Gamma^{FV}$, where $\Gamma^{FV}$ denotes the CPU time per leaf on the regular mesh, should be greater than one. 
Consequently, the  overhead per iteration and per leaf of the adaptive computation is defined by
\begin{equation}\label{eq:overhead}
 \mbox{\it overhead} = \frac{\Gamma^{a}}{\Gamma^{FV}}-1.
\end{equation}

\section{Numerical results}\label{numerical}

In the following, the results of the different MR and AMR  simulations are presented and compared with corresponding FV reference computations considering two Riemann problems, one in two and one in three space dimensions.

\subsection{2D Riemann problem: Lax-Liu configuration $\#6$}\label{sec:2dcomp}

\begin{figure}[t]\centering
\begin{tabular}{cc}
\parbox{5.5cm}{\includegraphics[height=4.5cm]{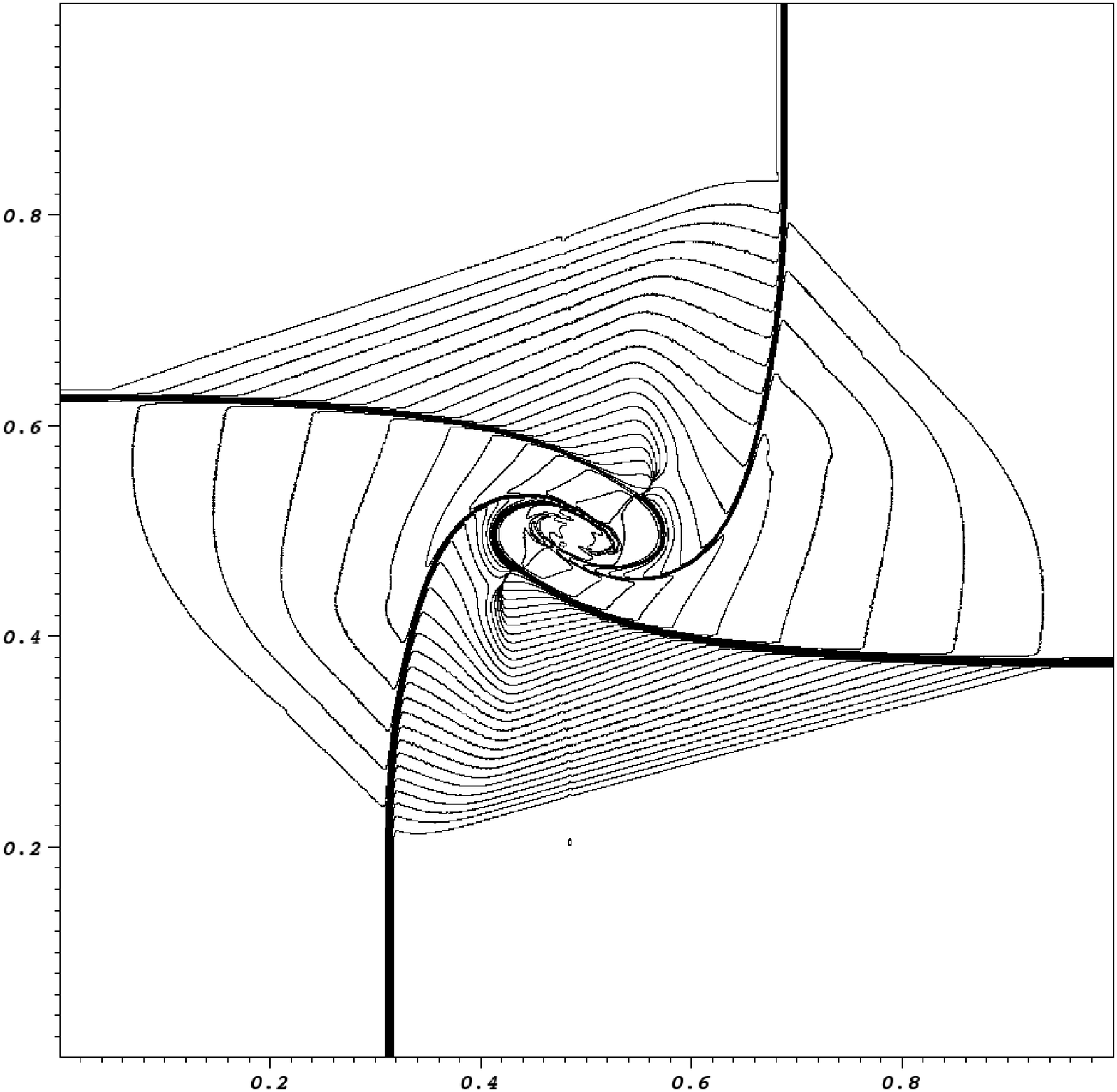}} &
\parbox{5.5cm}{\includegraphics[trim=2cm 2cm 2cm 2cm, clip=true, height=4.5cm,width=5.9cm]{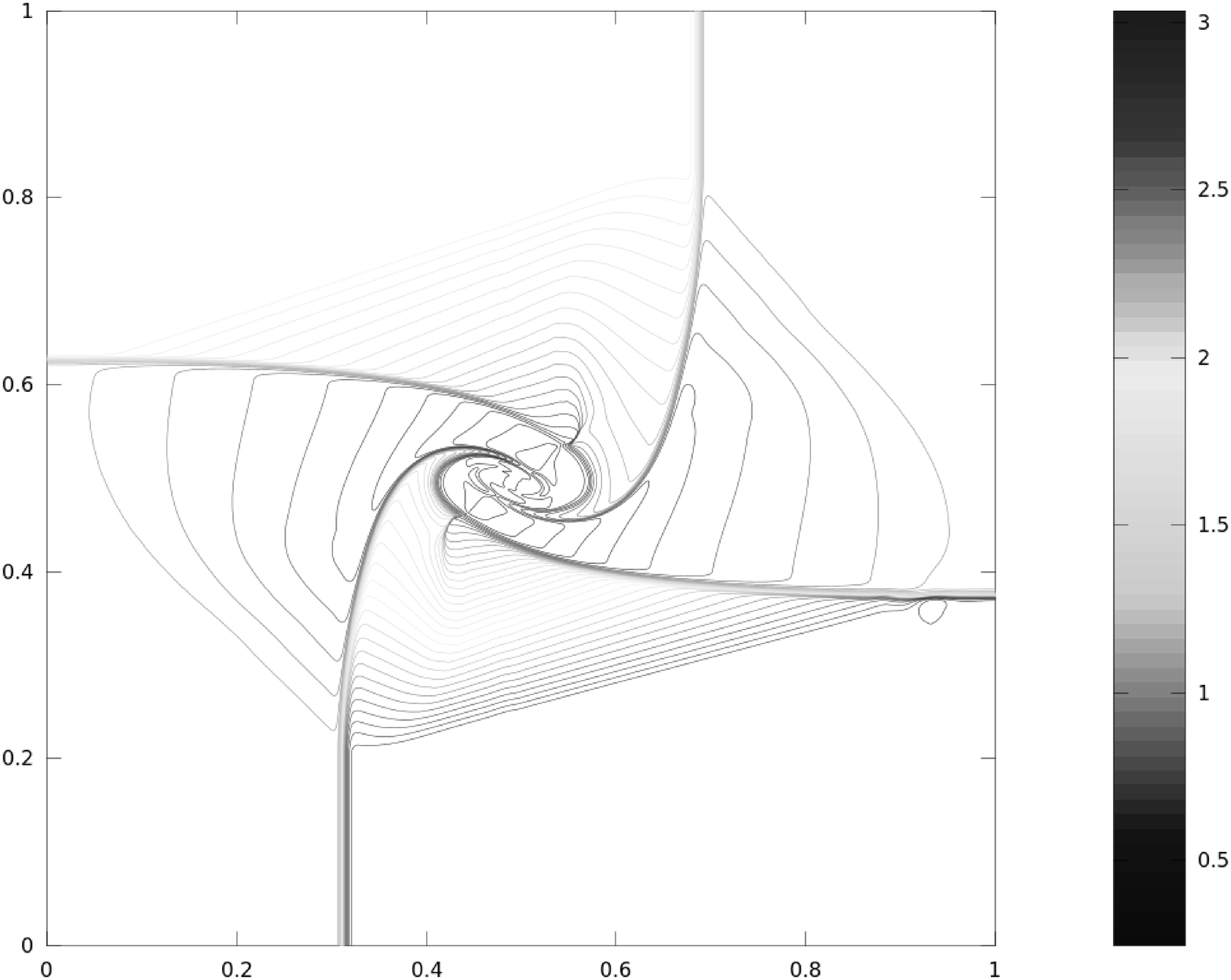}}
\end{tabular}
\caption{Left: reference solution -- isolines of density ($ 0\leq \rho \leq 3.2,$ with step $0.008$)  at the final time $t_e$ for the 2D Riemann problem Lax-Liu configuration \#6. 
The computation is performed on a $4096^2$ mesh using the FV volume scheme with AUSMDV flux. Right: MR solution, using $L=10$ and $\epsilon=0.0023$.}\label{pic:refausmdv}
\end{figure}

\begin{figure}[t]\centering
\begin{tabular}{cc}
\parbox{5.5cm}{\includegraphics[trim=2cm 0cm 2cm 0.5cm, clip=true, height=4.5cm]{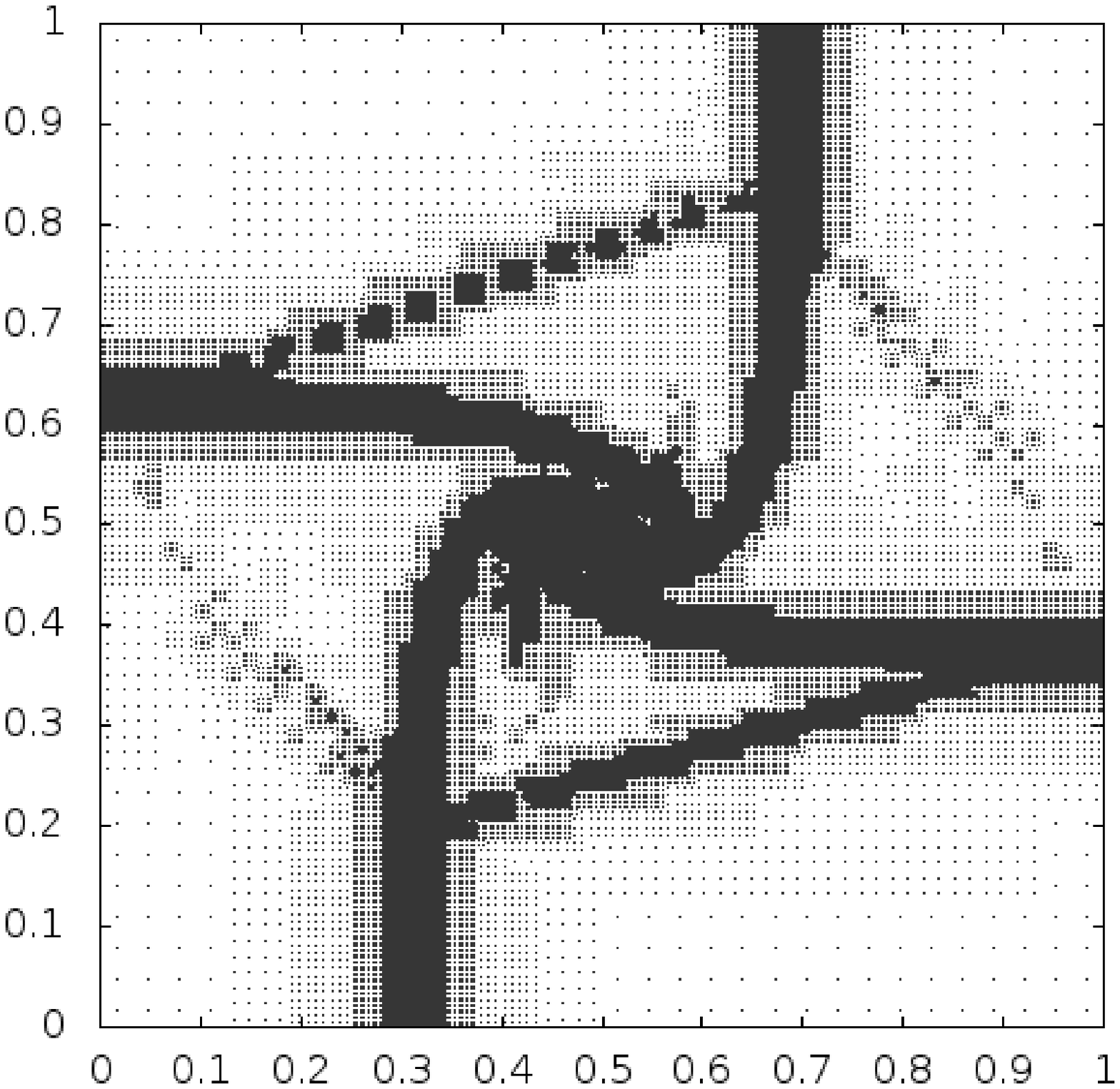}} & 
\parbox{5.5cm}{\vspace*{0.15cm}

\includegraphics[height=4.68cm]{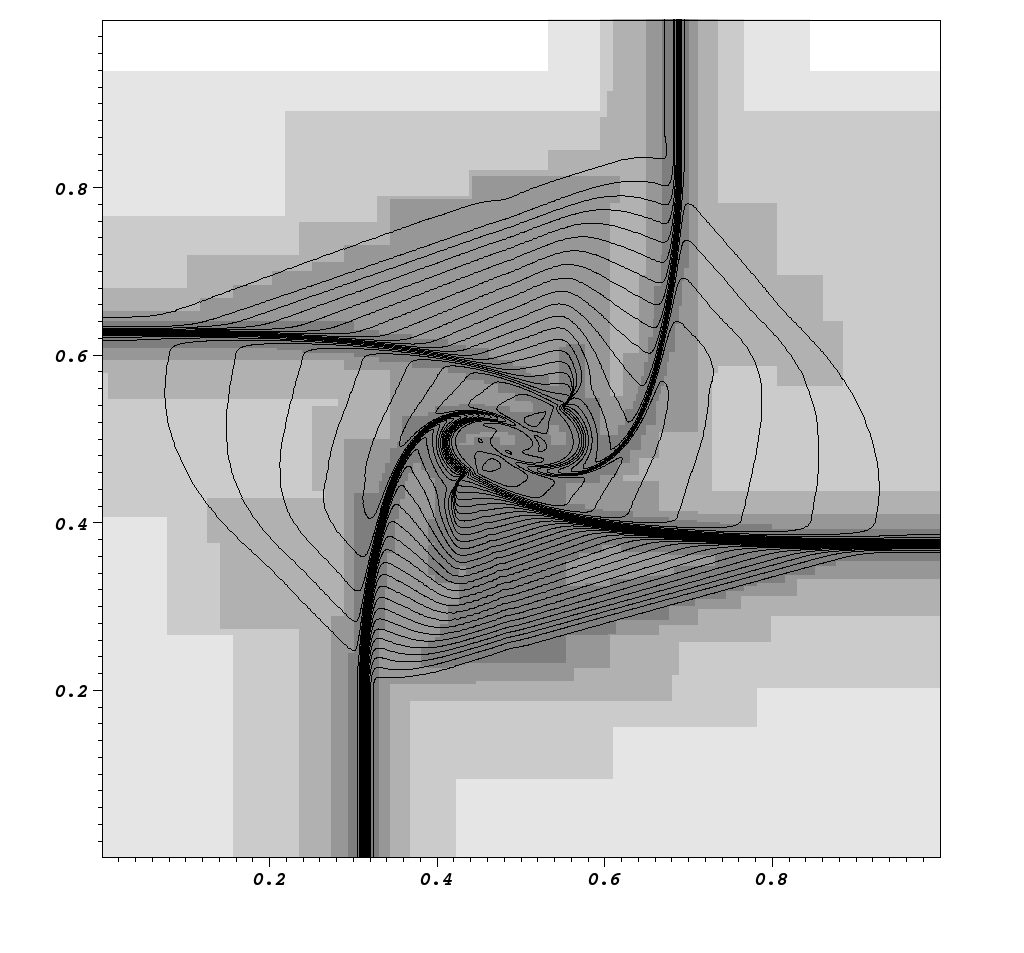}}   
\end{tabular}\vspace*{-0.2cm}
\caption{Adaptive solutions and meshes for the 2D Riemann problem Lax-Liu configuration \#6  at the final time $t_e$, using $L=10$:  
Left: MR mesh with $\epsilon=0.0023$. Right: AMR solution and mesh.}\label{pic:amrausmdv}
\end{figure}

First we consider a classical Riemann problem for gas dynamics proposed in \cite{LaxLiu:1998} known as
Lax-Liu benchmark \#$6$. This configuration is initially discussed in \cite{ZhangZheng:1990,SchulzRinneetal:1993}. \correction{As this is a 2D problem, $v_3\equiv 0$ holds true, the remaining vector components only depend on $x_1$ and $x_2$, and Eq.~(\ref{conservation}) reduces to $\partial_t {\bf q}+\partial_{x_1} f_1({\bf q})+\partial_{x_2} f_2({\bf q})=0$.}  
The computational domain is the square $\Omega=[0,1]^2$ with outflow boundary conditions. 
The domain is divided into four quadrants,  
where the initial data are set constant in each quadrant, according to the values  given in Table~\ref{tab:ICConf6}.  
The simulations are performed until the final time $t_e=0.25$. 
This benchmark only involves contact discontinuities and generates swirling motion in  the clockwise direction. 
\begin{table}[t]
\caption{Initial Values for the Lax-Liu configuration \#6. \label{tab:ICConf6}}
\centering\small
\begin{minipage}[c]{0.72\linewidth}
\centering\begin{tabular}{lcccc}
\hline \hline
\multicolumn{1}{c}{Variables} & \multicolumn{4}{c}{Domain position}\\ 
                & I    &   II &  III &  IV \\
\hline  
Density$(\rho)$ & 1.00 & 2.00 & 1.00 & 3.00 \\
Pressure $(p)$ & 1.00 & 1.00 & 1.00 & 1.00\\
Velocity component $(v_1)$   & 0.75 & 0.75 & -0.75 & -0.75 \\
Velocity component $(v_2)$   & -0.50 & 0.50 & 0.50 & -0.50 \\
\hline
\end{tabular}\vspace*{-0.3cm}
\end{minipage}
\begin{minipage}[c]{0.22\linewidth}
\centering
\def\pgfsysdriver{pgfsys-dvipdfm.def}
\begin{tikzpicture}[->,thick]
\draw[color=black] (0,0) --(2.3,0) node[right] {$x_1$} coordinate (x axis);
\draw[color=black] (0,0) --(0,2.3) node[above] {$x_2$} coordinate (y axis);
\draw[help lines]  (0,0) grid (2,2)
                   (0.5,0.5) node {III}
		   (1.5,0.5) node {IV}
		   (1.5,1.5) node {I}
		   (0.5,1.5) node {II}
                    (0,2.28)node{} ;
\end{tikzpicture}\vspace*{-0.3cm}
\end{minipage}
\end{table}
Figure~\ref{pic:refausmdv}, left shows the reference solution computed with the FV scheme of the AMR method on a regular mesh with $4096^2$ grid points which is used to 
evaluate the errors of the AMR and MR computations.\footnote{\correction{Computation of the reference solution on the uniform $4096^2$ grid used 
$8$ nodes with Intel-Xeon-$3.4$GHz dual-core processors of a typical GNU/Linux cluster and required $\sim 3.5\,\mathrm{h}$ wall time to complete $5120$ time steps. Note
that all benchmarked computations were run in serial.}}
Figures~\ref{pic:refausmdv}, right  and \ref{pic:amrausmdv}, right  illustrate, respectively,  the adaptive solutions obtained with MR  and AMR  at the final time $t_e$, using $L=10$ levels. 
The corresponding adaptive meshes are also shown in Figure~\ref{pic:amrausmdv}.
We observe that in the MR case the mesh is sparser and better adapted to the solution compared to the AMR case.

The underlying FV schemes of MR and AMR are benchmarked on uniform fine meshes with $N_C= N^2$ FV cells, where
$N=\left\{128,256,512,1024\right\}$, (corresponding to the scale levels $L=\left\{7,8,9, 10\right\}$). 
To reach the final instant $t_e$,  $N_I$ time steps are performed, where   $N_I= \left\{160, 320,640,1280\right\}$, respectively. 
Note that in all adaptive simulations $N_I$ time steps are performed as well. 
In both cases we observe convergence of the FV schemes towards the reference solution with a rate of about one, as shown in Table~\ref{tab:2Daccuracy}.
The FV scheme of MR yields slightly smaller errors with slightly higher convergence rates.
In the MR and MRLT computations we studied the influence of the wavelet threshold, and considered different values
 $\epsilon=0.0005, 0.0008, 0.0028,0.0025, 0.0023$, 0.0010, which were fixed for all levels.
For the MR and MRLT schemes, decreasing the threshold parameter $\epsilon$ has the effect of improving the accuracy to some extent, but mesh compression deteriorates for very small $\epsilon$.
In Table~\ref{tab:2Daccuracy} only results for $\epsilon=0.0023$ are presented.
Similarly to the FV case, we observe decaying errors for increasing $L$ with the same rates of about one.
This shows that the wavelet thresholding well preserves the order of the underlying FV scheme.

The AMR and AMRLT computations use a base mesh of 64 cells and a refinement factor of 2 at all levels. 
The full block-structured AMR algorithm is used here, including conservative correction at refinement boundaries and hierarchical time stepping. 
Refinement is based on the scaled gradient of the density with threshold $\epsilon_\rho=0.05$. 
Again, we observe, similar to the FV case, decaying errors for increasing $L$ with rates slightly below one.
The difference in the absolute errors between the  MR/MRLT and AMR/AMRLT simulations
is due to the use of slightly different FV schemes. In particularly, the FV errors of the AMR/AMRLT code are slightly larger, which is consistent
with the employed second order method, as discussed in Section~\ref{sec:nummethods}. However, the small perturbation rates of Table~\ref{tab:2Daccuracy} as well 
as the similar memory and mesh compression rates (differing only 
$2-4\%$) of Table~\ref{tab:2DMemory} indicate that both methods perform a comparable adaptation. 
It can be noted that the MR/MRLT compression rates decrease faster providing evidence for a more sophisticated adaptation criterion.

\begin{table}[t]
\caption{Accuracy of 2D computations with thresholds $\epsilon=0.0023 $ and $\epsilon_{\rho}=0.05$ in the MR and AMR cases, respectively.}\label{tab:2Daccuracy}\centering\small
\renewcommand{\arraystretch}{1.1}
\begin{tabular}{r||cc||cc|c||cc|c}
\hline\hline
\multirow{2}{*}{$L$} & \multicolumn{2}{|c||}{FV} &\multicolumn{2}{|c|}{MR} & Pert. & \multicolumn{2}{|c|}{MRLT} & Pert. \\
\cline{2-5} \cline{7-8} 
& $L_1(\rho)$ & Rate & $L_1(\rho)$   & Rate  & [\%] & $L_1(\rho)$   & Rate & [\%]  \\ \hline
$7$ &   $0.0390781$   & $ $     &  $0.0374855$  & $ $      & $4.03$ & $0.0376881$  &   $     $ &  $3.51$  \\
$8$ &   $0.0236107$   & $0.727$ &  $0.0225061$  & $0.736$  & $4.68$ & $0.0226217$  &   $0.736$ &  $4.19$  \\
$9$ &   $0.0128030$   & $0.883$ &  $0.0122309$  & $0.880$  & $4.47$ & $0.0123496$  &   $0.873$ &  $3.54$  \\
$10$ &  $0.0056620$   & $1.177$ &  $0.0056208$  & $1.122$  & $0.73$ & $0.0059899$  &   $1.044$ &  $5.79$  \\
\hline\hline  
\multirow{2}{*}{$L$} & \multicolumn{2}{|c||}{FV} &\multicolumn{2}{|c|}{AMR} & Pert. & \multicolumn{2}{|c|}{AMRLT} & Pert. \\
\cline{2-5} \cline{7-8} 
& $L_1(\rho)$ & Rate & $L_1(\rho)$   & Rate  & [\%] & $L_1(\rho)$   & Rate & [\%]  \\ \hline
$7$ &   $0.0458873 $  & $     $ &   $0.0436478 $ & $     $ & $4.88$ & $0.0456108$  & $     $ & $0.60$      \\
$8$ &   $0.0293788 $  & $0.643$ &   $0.0279279 $ & $0.644$ & $4.93$ & $0.0291994$  & $0.643$ & $0.61$    \\
$9$ &   $0.0174213 $  & $0.754$ &   $0.0165780 $ & $0.752$ & $4.84$ & $0.0173483$  & $0.751$ & $0.42$     \\
$10$ &   $0.0091630 $ & $0.927$ &   $0.0089999 $ & $0.881$ & $1.78$ & $0.0093855$  & $0.886$ & $2.43$     \\
\hline\hline  
\end{tabular}
\end{table}

\begin{table}[t]
\caption{Accumulated cell ($\cal C$) and leaf ($\cal L$) counts in $10^6$ elements of adaptive 2D computations and resulting 
memory and grid compression rates.}\label{tab:2DMemory}\centering \small
\renewcommand{\arraystretch}{1.1}
\begin{tabular}{r||cc|cc||cc|cc}\hline\hline
\multirow{2}{*}{$L$} & \multicolumn{4}{c||}{MR} & \multicolumn{4}{c}{MRLT}\\ \cline{2-9}
& $\sum_n {\cal C}^n$ & [\%] & $\sum_n {\cal L}^n$ & [\%] & $\sum_n {\cal C}^n$ & [\%] & $\sum_n {\cal L}^n$ & [\%] \\ \hline
$7$ & $2.16$ & $82.5$ & $1.46$ & $55.8$ & $2.12$ & $80.7$ & $1.43$ & $54.5$ \\
$8$ & $11.2$ & $53.5$ & $7.35$ & $35.1$ & $11.0$ & $52.3$ & $7.19$ & $34.3$ \\
$9$ & $50.3$ & $30.0$ & $32.4$ & $19.3$ & $49.9$ & $29.8$ & $32.2$ & $19.2$ \\
$10$&  $210$ & $15.7$ &  $135$ & $10.1$ &  $207$ & $15.4$ &  $133$ & $9.91$ \\
\hline\hline  
\multirow{2}{*}{$L$} & \multicolumn{4}{c||}{AMR} & \multicolumn{4}{c}{AMRLT}\\ \cline{2-9}
& $\sum_n {\cal C}^n$ & [\%] & $\sum_n {\cal L}^n$ & [\%] & $\sum_n {\cal C}^n$ & [\%] & $\sum_n {\cal L}^n$ & [\%] \\ \hline
$7$ & $1.84$ & $70.3$ & $1.39$ & $52.9$ & $1.88$ & $71.8$ & $1.28$ & $49.0$ \\
$8$ & $9.79$ & $46.7$ & $7.35$ & $35.0$ & $10.3$ & $48.9$ & $6.80$ & $32.4$ \\
$9$ & $50.1$ & $29.9$ & $37.6$ & $22.4$ & $51.4$ & $30.6$ & $33.3$ & $19.9$ \\
$10$&  $239$ & $17.8$ &  $180$ & $13.4$ &  $242$ & $18.0$ &  $154$ & $11.4$ \\
\hline\hline  
\end{tabular}
\end{table}

Table~\ref{tab:cpu} presents the absolute computing time of both codes and the corresponding CPU time compression rates. Comparing the absolute CPU times already for
the FV unigrid cases, it is important to point out that the two benchmarked implementations have vastly different absolute computational performance. 
When running in FV mode, the MR/MRLT code is about a factor of $9.5$ slower than the AMR/AMRLT code and an accordingly larger absolute
computing time is consequently also measured in the adaptive simulations. In order to allow nevertheless a comparison, and to assess both mathematical approaches independent of mere implementation aspects \correction{as well as minor differences in the numerical discretizations applied}, we employ the CPU time compression rate. 
For the highest resolved MR/MRLT simulations the CPU time
compression is slightly better than for AMRLT, which is consistent with the better mesh compression rates seen in Table~\ref{tab:2DMemory}, and thus significantly
better than for AMR. Contrary to the AMRLT computations, the present MRLT scheme does not improve the results significantly with respect to MR.
This issue of local time stepping in MR is also discussed in \cite{DGRS:2008}.
MRLT is most beneficial for unbalanced trees, \textit{i.e.}, localization of small scale features of the solution, expensive flux evaluation, and larger number of 
well-localized active scales. Therefore, for few active cells on fine scales (\textit{e.g.}, point singularities) the speed-up becomes larger compared to the MR scheme with global time stepping.
Similar results are also found by  M\"uller and Stiriba \cite{MullerStiriba:2007} for their MRLT computations for space dimensions larger than one. 


As predicted by an analytic cost estimate in \cite{DGRS:2008}, we found that the actual speed-up of local time stepping depends on the distribution 
of the active cells. If the majority of the cells is active on fine scales, the MRLT scheme is less efficient with respect to the MR scheme, whereas 
for few active cells on fine scales the speed-up becomes larger. 
In these cases, performing one large time step at a coarse level instead of several time steps on fine scale
cells becomes indeed more efficient.

\correction{When analyzing overhead rates according to Eq.~(\ref{eq:overhead}), one notices that the relative mesh adaptation overhead rises in all approaches when
increasing the number of levels $L$. The respective variants without local time stepping exhibit a steeper increase than the LT methods as coarser grid cells are 
updated significantly more often. For the highest refined computations with $L=10$, the overhead is $158\%$ for AMR and $77.0\%$ for AMRLT. When the overhead is calculated 
based on the performance of the MR code run in FV unigrid mode, the corresponding overhead rates for MR and MRLT are only $72.3\%$ and $67.5\%$, respectively. 
Measured in relative terms, the MR/MRLT algorithms are more efficient than AMR/AMRLT, which is due to the unnecessary update of respective coarser cells 
independent of higher level coverage. However, if the performance numbers from the AMR code run in FV unigrid mode would be used for the sake of  
direct comparison, the overhead rates for MR and MRLT would jump to $1560\%$ and $1513\%$, respectively. This underscores that exactly identical 
numerical discretizations and implementations in the same programming language of comparable quality and level of optimization are a prerequisite for
meaningful head to head benchmarking of adaptive solution methods. Yet, due to the complexity of the involved codes, these requirements can generally 
never be met in practice, which explains the lack of published quantitative comparisons in this field and provides evidence for the benefit of our
approach using primarily relative performance measures.}

\begin{table}[t]
\caption{Computing times and CPU time compression rates for 2D computations.}\label{tab:cpu}\centering\small
\renewcommand{\arraystretch}{1.1}
\begin{tabular}{r||c|cc|cc||c|cc|cc}
\hline\hline
\multirow{2}{*}{$L$} & FV & \multicolumn{2}{c|}{MR} & \multicolumn{2}{c||}{MRLT}      & FV & \multicolumn{2}{c|}{AMR}  & 
\multicolumn{2}{c}{AMRLT}\\ \cline{3-6}\cline{8-11}
&[s] &[s] & [\%]&[s] & [\%]         &[s] &[s] & [\%]  & [s] & [\%]\\
 \cline{2-7}
\hline
$7$ & $27.22$ & $24.42$ & $89.7$ & $23.61$ & $86.7$ &  $2.901$ & $3.864$ &  $133$ & $2.489$ & $86.2$  \\
$8$ & $212.4$ & $124.4$ & $58.6$ & $118.7$ & $55.9$ &  $23.36$ & $19.76$ & $85.0$ & $12.64$ & $54.1$  \\
$9$ &  $1683$ & $555.0$ & $33.0$ & $538.0$ & $32.0$ &  $181.6$ & $99.10$ & $54.6$ & $60.90$ & $33.5$  \\
$10$& $13727$ &  $2380$ & $17.3$ &  $2278$ & $16.6$ &  $1420 $ & $489.3$ & $34.5$ & $287.4$ & $20.2$  \\
\hline\hline  
\end{tabular}
\end{table}

\begin{table}[t]
\caption{Breakdown of CPU time in \% spent in main task groups for 2D computations at $L=10$.}
\label{tab:breakdown}
\centering\small
\renewcommand{\arraystretch}{1.1}
\begin{tabular}{l|rrr||l|rrr}\hline\hline
\multirow{2}{*}{Task group} & \multicolumn{3}{c||}{MR} & \multirow{2}{*}{Task group} & \multicolumn{3}{c}{AMR} \\ \cline{2-4}\cline{6-8}
& \multicolumn{1}{c}{FV} & \multicolumn{1}{c}{NLT} & \multicolumn{1}{c||}{LT} & & \multicolumn{1}{c}{FV} & 
\multicolumn{1}{c}{NLT} & \multicolumn{1}{c}{LT} \\ \hline 
Numerics           & 36.18 & 26.30 & 25.30 & Numerics        & 93.73 & 64.44 & 71.55 \\ \cline{5-8}
Temp. data     & 33.06 & 16.67 & 16.62 & AMR data org.   &    -- &  6.72 &  4.16 \\ 
Boundary cond.     &  5.52 &  2.60 &  2.55 & Clustering      &    -- & 10.84 &  6.86 \\ \cline{1-4}
MR tree org.  &    -- & 40.63 & 40.28 & Flagging        &    -- &  2.32 &  1.32 \\ 
Level transfer     &    -- &  1.93 &  2.03 & Lever transfer  &    -- &  4.92 &  4.08 \\ \hline
Memory (c-lib)     & 21.79 &  6.45 &  6.23 & Memory (c-lib)  &  4.87 &  4.47 &  5.11 \\             
Unassigned         &  3.45 &  5.42 &  6.45 & Unassigned      &  1.40 &  6.29 &  6.92 \\ \hline \hline
\end{tabular}
\end{table}

\correction{In order to better understand nevertheless the absolute performance of the used computer codes, a breakdown of the CPU time of}
the highest resolved adaptive computations at $L=10$ is provided in Table~\ref{tab:breakdown}. This benchmark analysis has been
obtained by using the {\tt perf} tool. The tabulated items in the first group list major tasks already present in the respective FV implementations.
The second group lists algorithmic tasks that have to be carried out in addition in the adaptive programs as they could
be distinguished and unambiguously classified from within the performance analysis tool. The last group lists the expense of system 
calls for dynamic memory management and the accumulation of operations which have too insignificant cost for detailed classification. 
It can be seen from Table~\ref{tab:breakdown} that temporary data creation, handling and deletion are a considerable portion of the FV-MR code. This
code is a cell-oriented implementation of a Cartesian unstructured grid, while the FV method in the AMR code
is implemented for a generic structured block and allocates all necessary temporary data in advance. 
In the adaptive case, creation and organization of the hierarchical data
becomes a visible portion of the computing time for both approaches. It can be seen that particularly orchestration of the unstructured quad-tree in the MR/MRLT approach,
cf. Fig.~\ref{fig:tree}, is rather expensive. In principle, organizing the more general mesh refinement tree in the AMR/AMRLT case is even more involved; however, 
thanks to clustering individual cells into blocks of considerable size the number of leaves of this general tree, cf. Fig.~\ref{pic:hier}, is magnitudes smaller
and the overhead accordingly reduced. A non-negligible portion of the AMR approach is the {\em clustering} operation, which generates rectangular blocks from 
ensembles of individual cells. It is obvious that in the variant \correction{without local time step} refinement (AMR) the adaptation overhead is increased, since in the runs 
reported here complete mesh adaptation is allowed in every time step. 


\subsection{3D ellipsoidally expanding shock wave}\label{sec:3dcomp}
For the next test problem, we consider the expansion of an ellipsoidal shock wave in three space dimensions. 
The 3D Euler equations are solved in the computational domain $\Omega=[-2,2]^3$ until the final simulation time $t_e=0.28$. Outflow boundary conditions
are applied at all sides of the domain. 
The initial ellipsoid is specified by 
\begin{equation}
 r=\sqrt{\left(\frac{x_{1_r}}{a}\right)^2+\left(\frac{x_{2_r}}{b}\right)^2+\left(\frac{x_{3_r}}{c}\right)^2},
 \end{equation}
where
$
x_{1_r}= x_1\cos(\theta)-x_2\sin(\theta), \;\;
x_{2_r}=\left(x_1\sin(\theta)+x_2\cos(\theta)\right)\cos(\phi)-x_{3}\sin(\phi), 
\newline x_{3_r}=\left(x_1\sin(\theta)+x_2\cos(\theta)\right)\sin(\phi)+x_3\cos(\phi),
$
with stretching and rotational parameters $r_c=\frac{3}{5}$, $a=\frac{1}{3}$, $b=1$, $c=3$, $\theta=\frac{\pi}{3}$, and $\phi=\frac{\pi}{4}$. Initial conditions
in density $\rho$ and energy density $\rho e$ are set as 
\begin{equation}
 \rho = \left\{ \begin{array}{ll} 0.125 \,, & r<r_c\;, \\ 1\,, &  r\ge r_c \;, \end{array} \right. \qquad 
\rho e = \left\{ \begin{array}{ll} 0.25 \,, & r<r_c\;, \\ 2.5\,, &  r\ge r_c \;, \end{array} \right.  
\end{equation}
while the velocity vector is initially zero, \textit{i.e.}, $v_{1}=v_{2}=v_{3}=0$ everywhere.
For convergence analysis, we consider a reference solution computed with AUSMDV numerical flux and second order MUSCL-Hancock reconstruction 
with Minmod-limiter in the primitive variables $\rho$, $v_n$, $p$. The resolution is $1024^3$ cells on a uniform mesh, where automatic time step adjustment based on 
$\mathrm{CFL}\approx 0.4$ is used.\footnote{\correction{While all benchmarked 3D simulations where run in serial, only the computation of this reference solution on the uniform $1024^3$ grid 
was performed on $512$ processors on an IBM-BG/P used in SMP mode, which required $\sim 9.2\,\mathrm{h}$ wall time to complete $460$ time steps.}} 

\begin{figure}[t] 
\centering
\begin{tabular}{ccc}
\includegraphics[width=3.15cm]{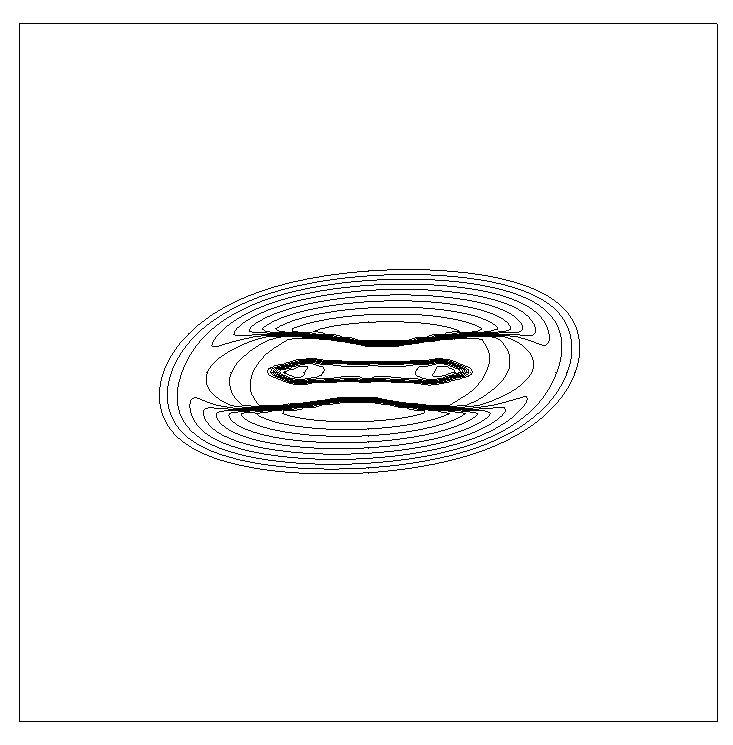}\hspace*{-0.5cm} &
\includegraphics[width=3.15cm]{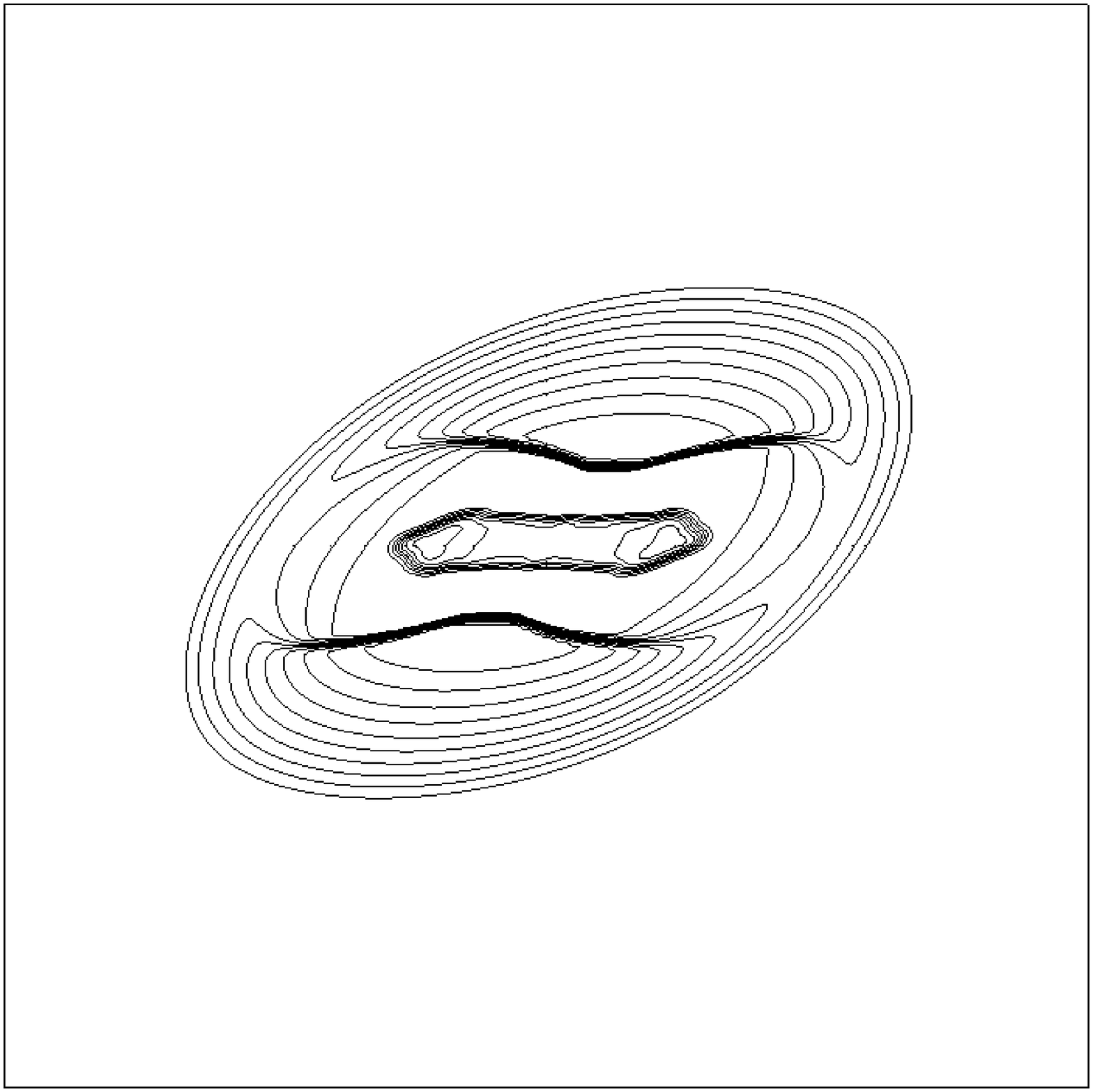}\hspace*{-0.5cm} &
\includegraphics[width=3.15cm]{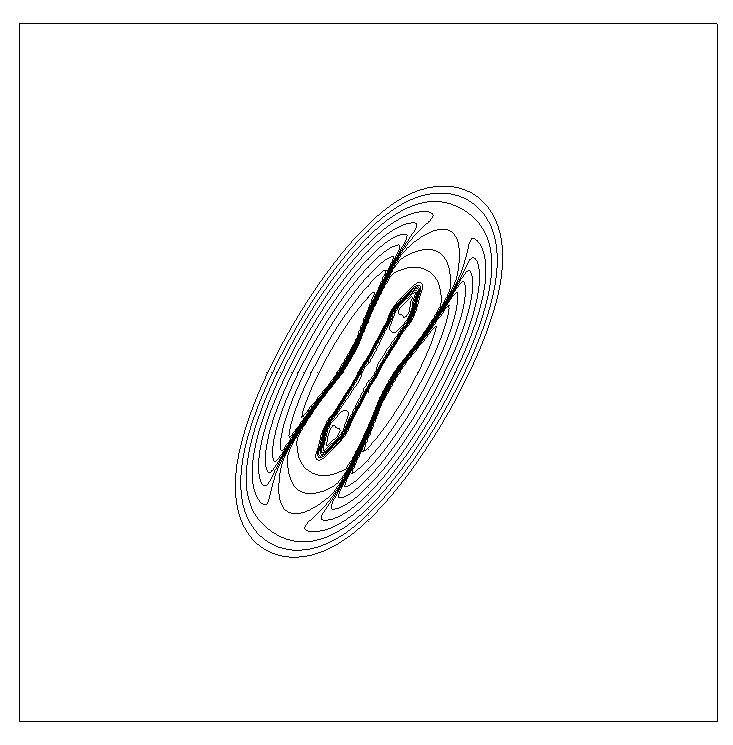} \\
\includegraphics[width=3.825cm,height=3.6cm]{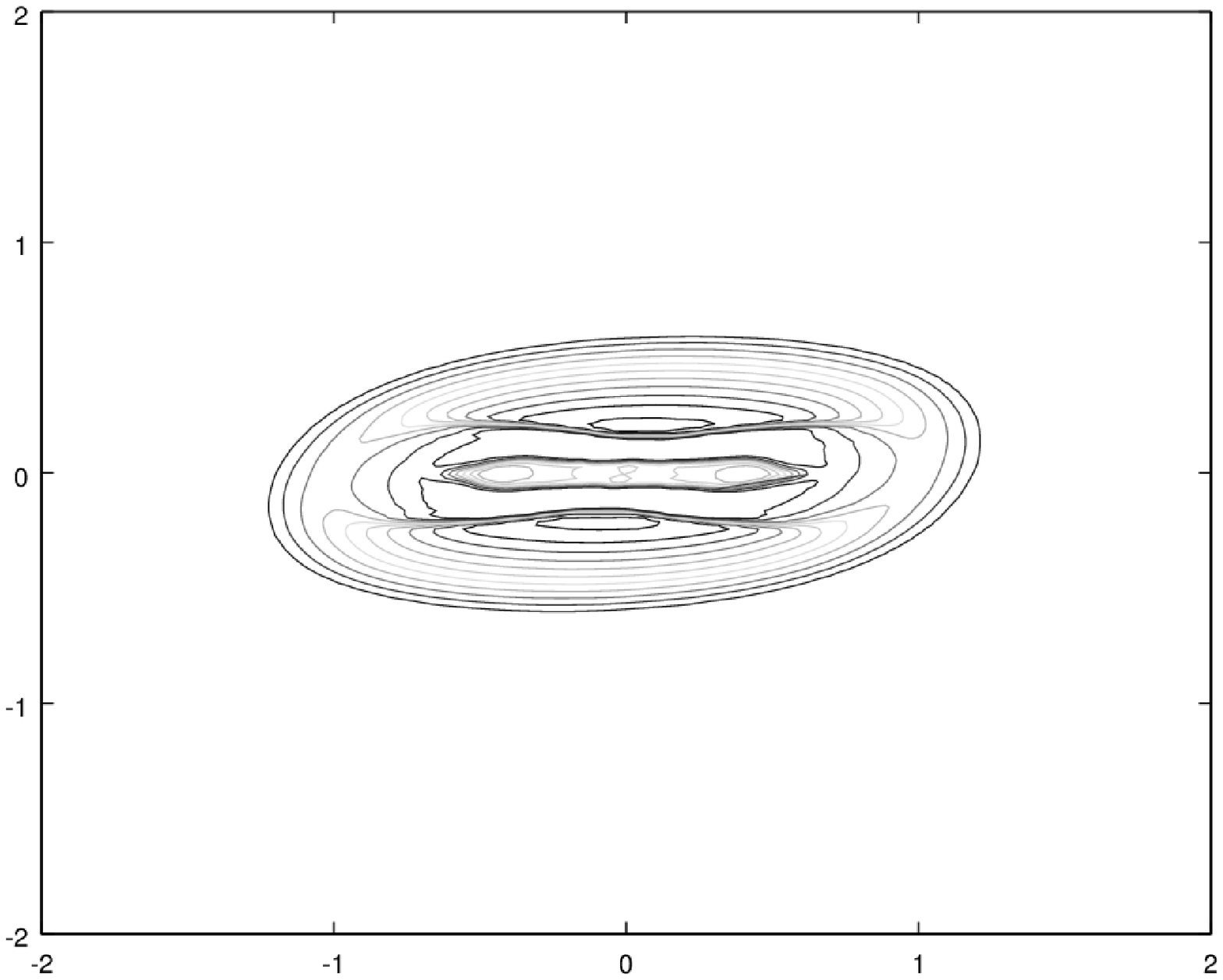}\hspace*{-0.4cm} &
\includegraphics[width=3.825cm,height=3.6cm]{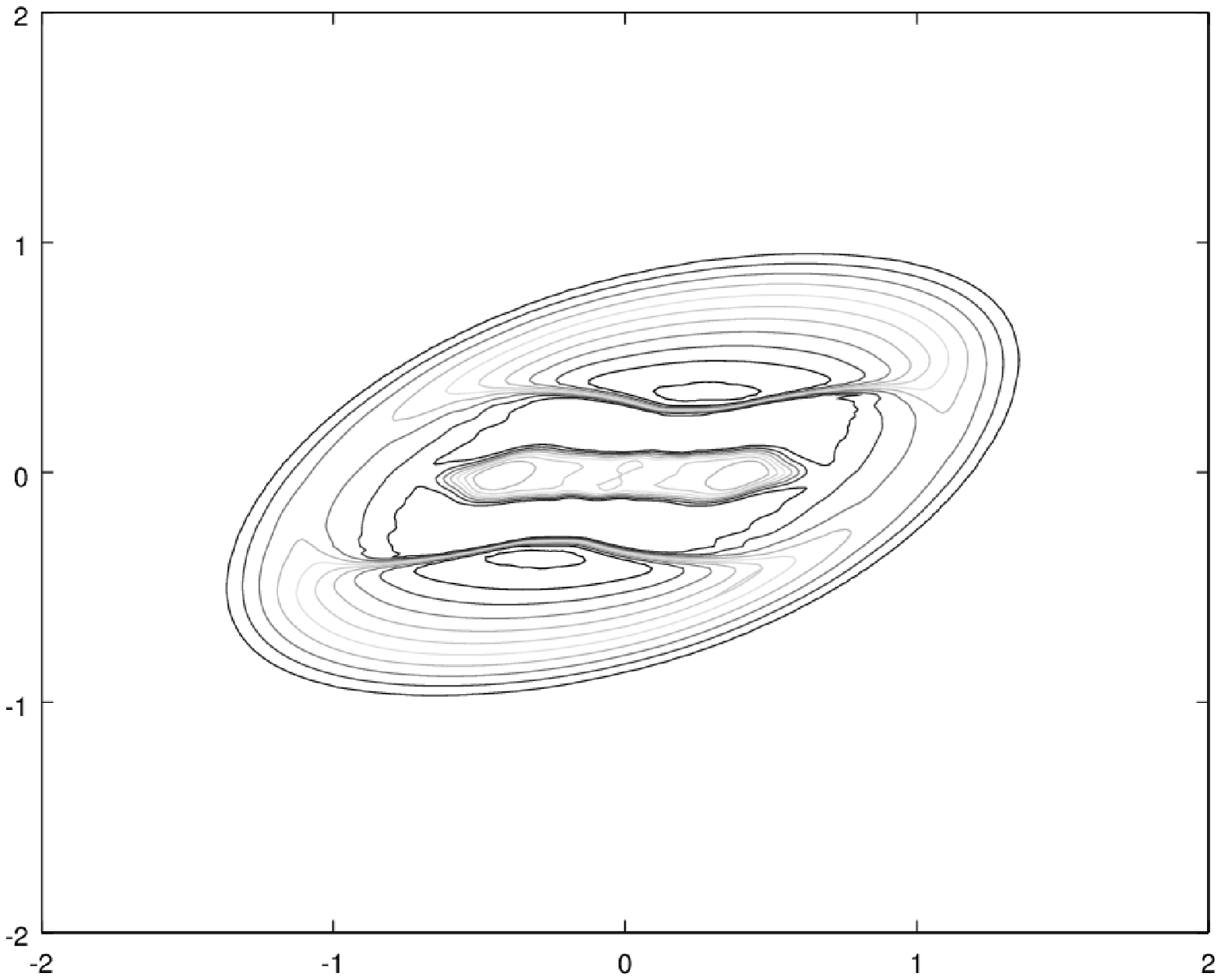}\hspace*{-0.4cm} &
\includegraphics[width=3.825cm,height=3.6cm]{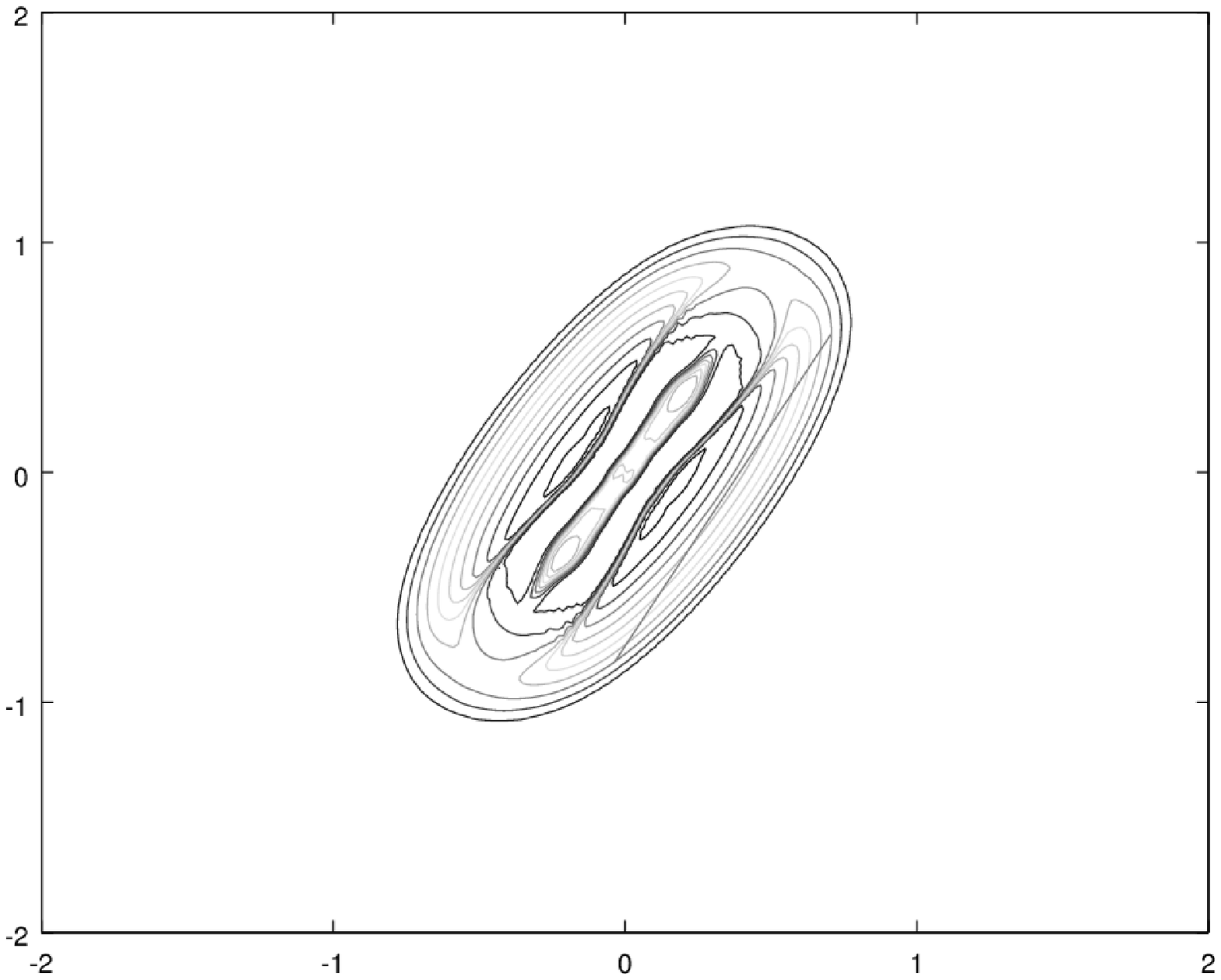} 
\end{tabular}\vspace*{-0.2cm}
\caption{3D ellipsoidally expanding shock-wave at time $t_e$.
Upper row: Isolines of 2D cuts of density for the
  FV reference solution computed on $1024^3$ mesh, down-sampled to $256^3$, corresponding to $L=8$.
Lower row: Isolines of 2D cuts of density for the MR computation with $L=8$.}
\label{fig:3dref}
\end{figure}

\begin{figure}[t]
\centering
\begin{tabular}{ccc}
\includegraphics[width=3.6cm]{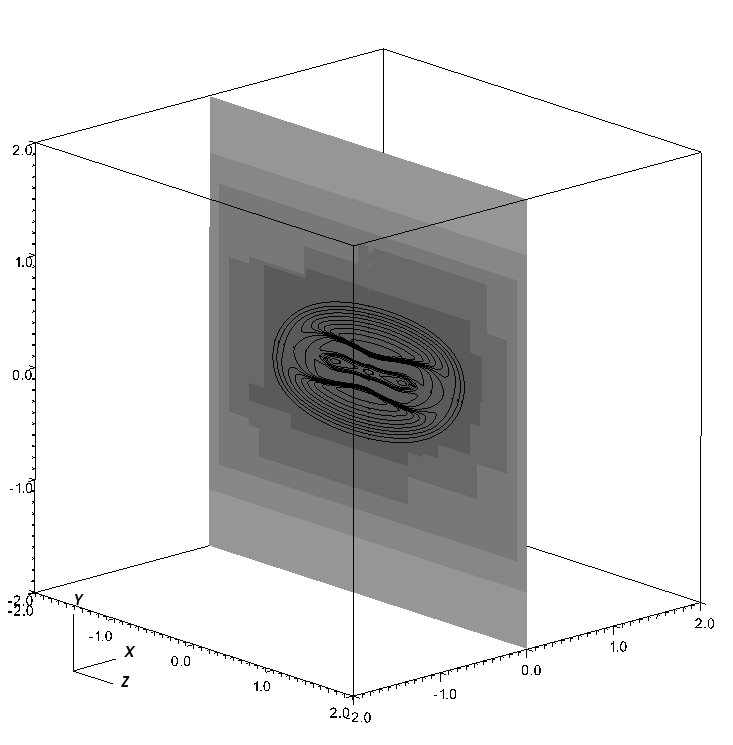}\hspace*{-0.1cm} & 
\includegraphics[width=3.6cm]{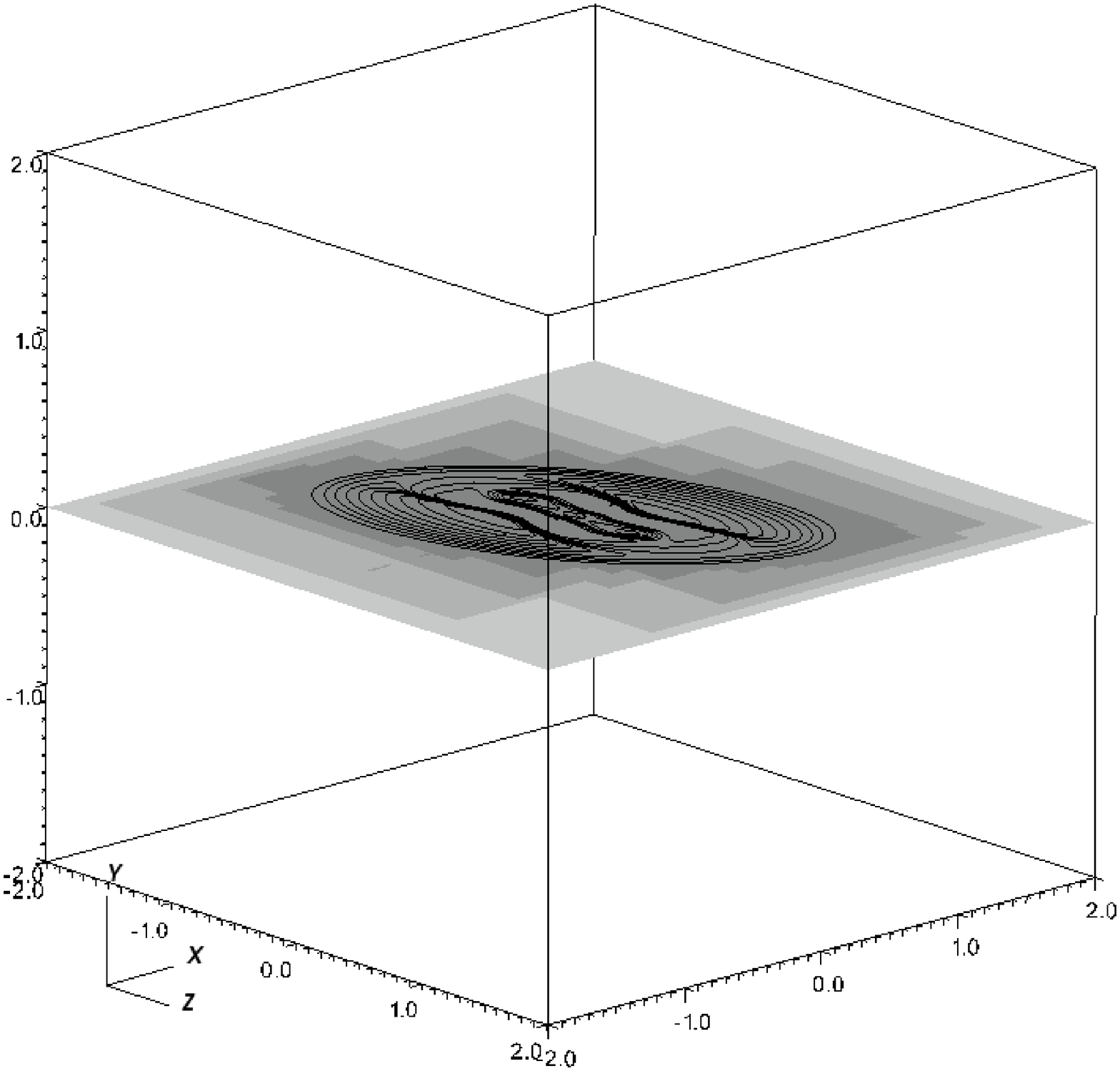}\hspace*{-0.1cm} &
\includegraphics[width=3.6cm]{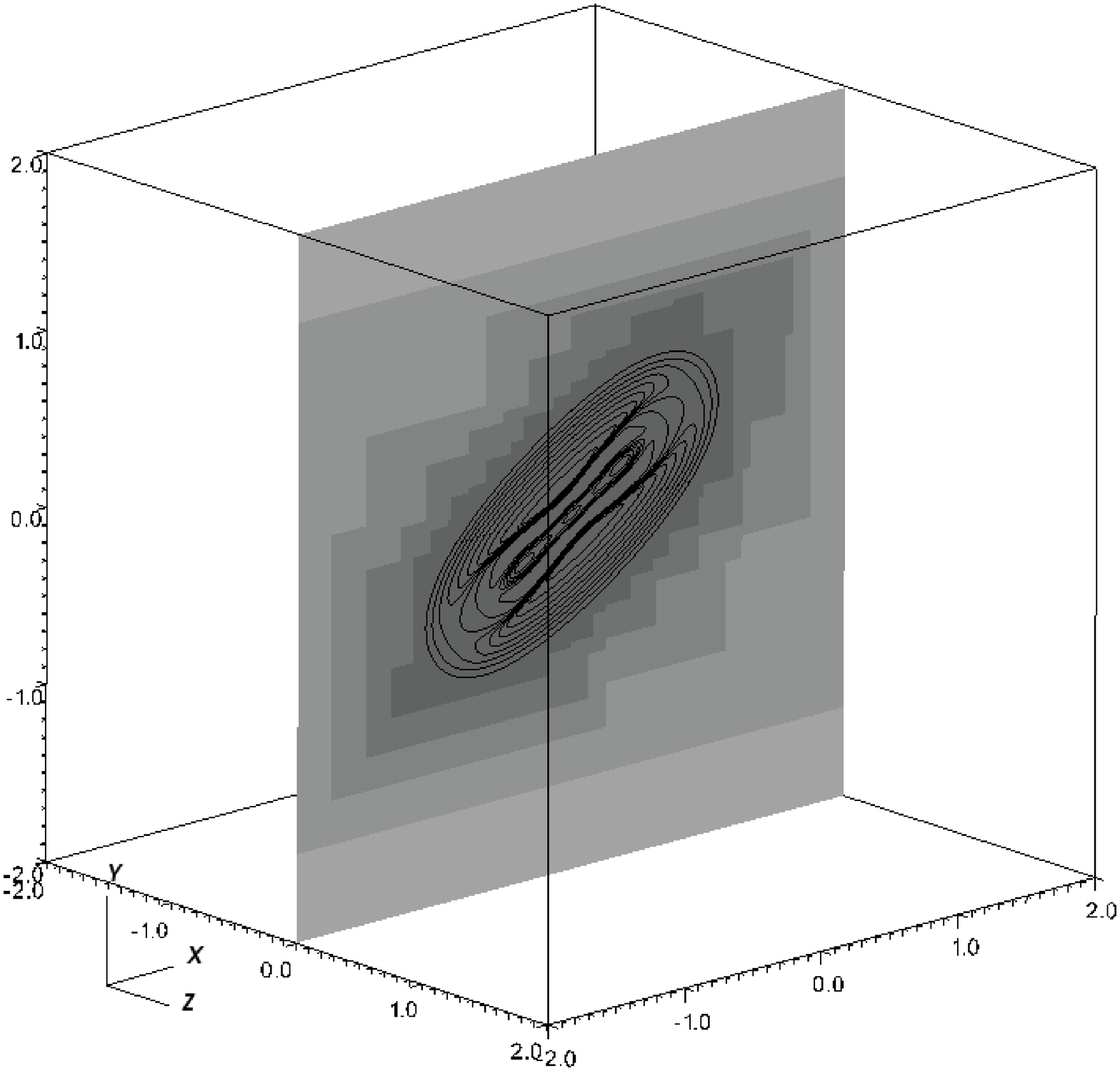}\hspace*{-0.1cm} \\
\includegraphics[width=3.42cm]{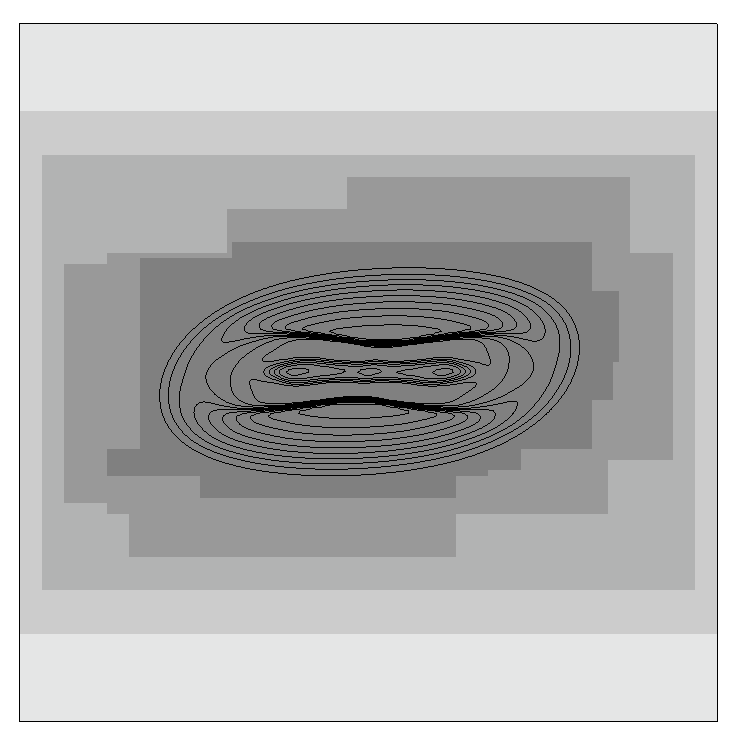}\hspace*{-0.1cm} & 
\includegraphics[width=3.42cm]{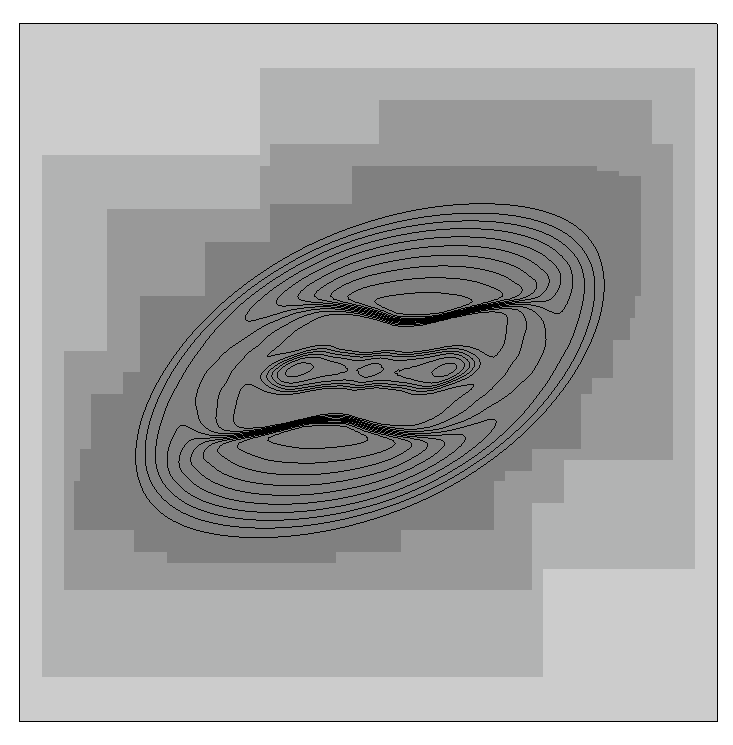}\hspace*{-0.1cm} & 
\includegraphics[width=3.42cm]{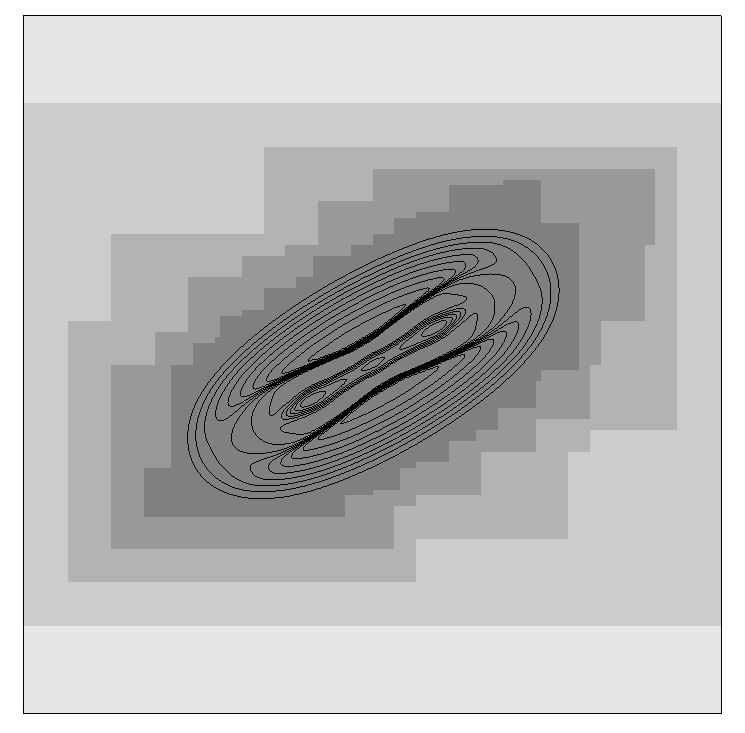}\hspace*{-0.1cm} \\ 
\hspace*{-0.3cm}\includegraphics[trim=2.1cm 0cm 2.1cm 0cm, clip=true, width=3.69cm]{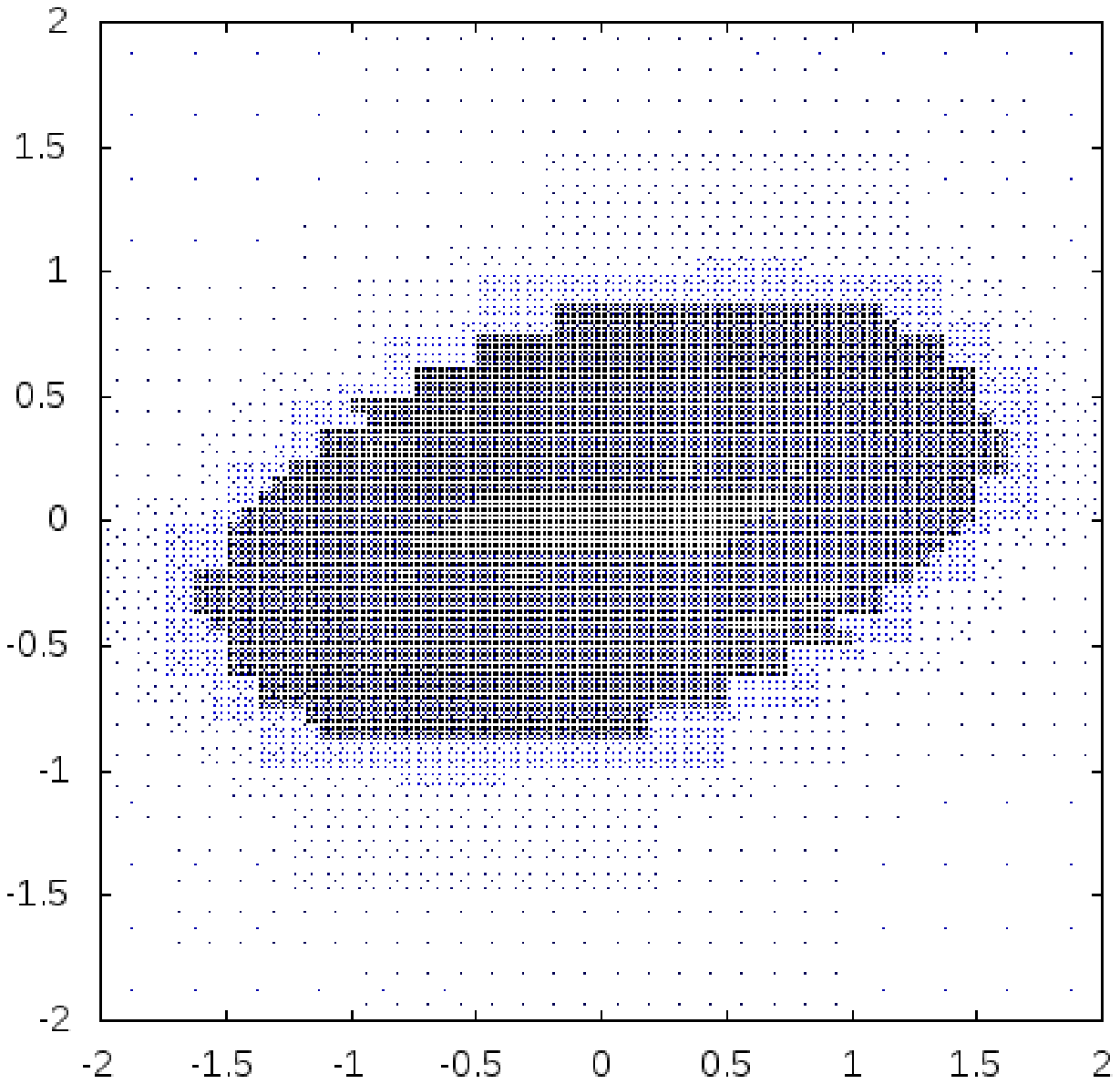}\hspace*{-0.1cm}  & 
\hspace*{-0.3cm}\includegraphics[trim=2.1cm 0cm 2.1cm 0cm, clip=true, width=3.69cm]{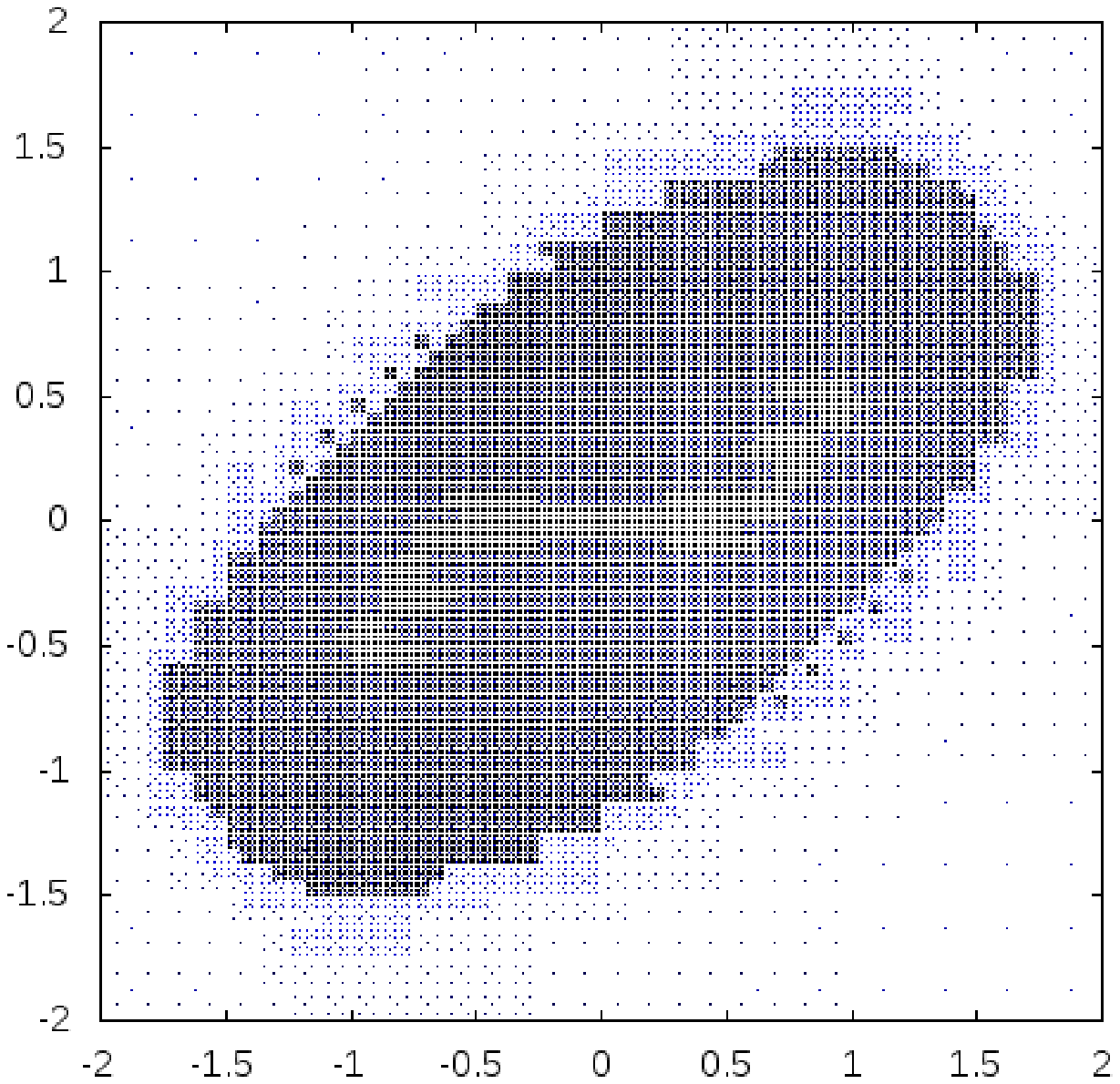}\hspace*{-0.1cm}  & 
\hspace*{-0.3cm}\includegraphics[trim=2.1cm 0cm 2.1cm 0cm, clip=true, width=3.69cm]{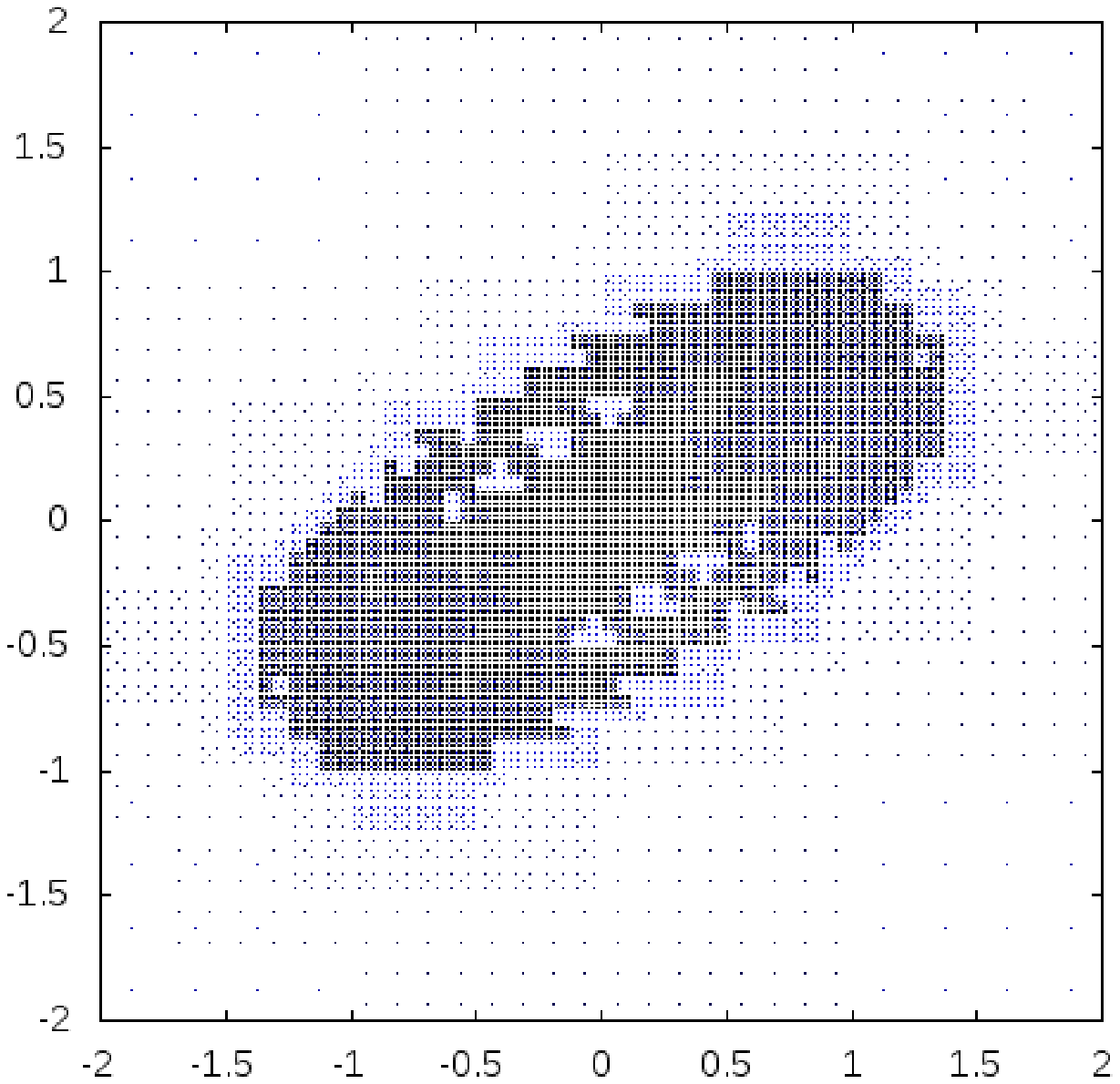}\hspace*{-0.1cm}  \\ 
\end{tabular}
\caption{Adaptation for 3D shock-wave at time $t_e$.
Upper row: AMRLT computation with $L=8$. 3D visualization of 2D cuts of density in the different planes 
( $y-z$ plane  at $x=0$,  $x-z$ plane  at $y=0$,  $x-y$ plane  at $z=0$),
%
Middle row: Isolines of 2D cuts of density for the AMRLT computation with $L=8$, superimposing the refinement levels in gray.
Lower row: Subset of the adaptive MR mesh  with $L=8$  projected onto the planes shown in middle row, 
considering all cells with positions in the interval $[-0.35,0.35]$ in the perpendicular direction of the respective plane.}
\label{fig:3dall}
\end{figure}

%
The AMR and AMRLT computations use a base mesh of $8^3$ cells and a refinement factor $2$ at all levels. 
As in the 2D example, the entire block-structured AMR algorithm is applied, including a conservative correction at the refinement boundaries and hierarchical time stepping. 
The refinement is based on scaled gradients of density and pressure with the thresholds 
$\epsilon_\rho=\epsilon_p=0.05$. The used mesh generation efficiency is $\eta_{tol}=80\,\%$, meaning that each patch can contain up to $20\,\%$ of cells not flagged for refinement. 
For the MR and MRLT computations, the threshold parameter $\epsilon=0.013$ is applied at all resolution levels.

%
The Figs.~\ref{fig:3dref} and \ref{fig:3dall} assemble 2D cuts normal to the coordinate directions at the final time  $t_e$. 
The left image of Fig.~\ref{fig:3dall} shows the plane in the origin in the view direction $(1,0,0)$, the middle one displays the plane in the 
origin in the view direction $(0,-1,0)$, and the right image shows the plane through the origin in the view direction $(0,0,-1)$.
In Fig.~\ref{fig:3dref}, we compare the FV reference computation obtained with $L=10$ and down-sampled to $L=8$ with the MR computations for $L=8$.
Shown are only isolines of density ($\rho$ between $0.05$ and $0.95$ at intervals $0.05$). Figure~\ref{fig:3dall} displays the 2D cuts and 
the corresponding adaptive meshes.
Inspecting these graphics, we find that both the MR and the AMRLT adaptive computation agree well with the down-sampled FV reference solution. 
As expected, both adaptive grids concentrate cells in regions of steep gradients present in the solution. 
We also observe that the MR mesh is better adapted to the solution than the block-structured AMR mesh, \textit{i.e.}, the MR mesh is sparser.
Similar to the previous section, we perform a detailed analysis of the computations. 
The underlying FV schemes of the adaptive methods are benchmarked on uniform meshes first with $N_C=N^3$ cells, where $N=\left\{32,64,128,256\right\}$, 
which corresponds to the levels $L=\left\{5,6,7,8\right\}$. To reach the final time instant $t_e$, $N_I$ time steps are performed, 
where $N_I= \left\{8, 32, 64, 128\right\}$, respectively. The number of time steps $N_I$ is the same in all adaptive simulations.

\begin{table}[t]
\caption{Accuracy of the 3D computations for the two FV schemes, MR and MRLT, and AMR and AMRLT using different levels of refinement.}\label{tab:3Accuracy}\centering\small
\renewcommand{\arraystretch}{1.1}
\begin{tabular}{r||cc||cc|c||cc|c}
\hline\hline
\multirow{2}{*}{$L$} & \multicolumn{2}{|c||}{FV} &\multicolumn{2}{|c|}{MR} & Pert. & \multicolumn{2}{|c|}{MRLT} & Pert. \\
\cline{2-5} \cline{7-8} 
& $L_1(\rho)$ & Rate & $L_1(\rho)$   & Rate  & [\%] & $L_1(\rho)$   & Rate & [\%]  \\ \hline
$5$ &   $0.386740$  &         &  $0.387222$ &         & $0.12$ & $0.387230$ &         & $0.13$  \\
$6$ &   $0.227130$  & $0.768$ &  $0.227270$ & $0.769$ & $0.06$ & $0.227278$ & $0.769$ & $0.07$  \\
$7$ &   $0.121120$  & $0.907$ &  $0.121157$ & $0.908$ & $0.03$ & $0.121248$ & $0.906$ & $0.11$  \\
$8$ &   $0.053677$  & $1.174$ &  $0.054001$ & $1.166$ & $0.60$ & $0.055429$ & $1.129$ & $3.26$  \\
\hline\hline  
\multirow{2}{*}{$L$} & \multicolumn{2}{|c||}{FV} &\multicolumn{2}{|c|}{AMR} & Pert. & \multicolumn{2}{|c|}{AMRLT} & Pert. \\
\cline{2-5} \cline{7-8} 
& $L_1(\rho)$ & Rate & $L_1(\rho)$   & Rate  & [\%] & $L_1(\rho)$   & Rate & [\%]  \\ \hline
$5$ &   $0.472266$  &         & $0.472295$ &         & $0.00$ & $0.472313$ &         & $0.00$  \\
$6$ &   $0.309526$  & $0.610$ & $0.309659$ & $0.609$ & $0.04$ & $0.309662$ & $0.609$ & $0.04$  \\
$7$ &   $0.184180$  & $0.749$ & $0.184503$ & $0.747$ & $0.18$ & $0.184430$ & $0.748$ & $0.14$  \\
$8$ &   $0.086319$  & $1.093$ & $0.087123$ & $1.083$ & $0.93$ & $0.086900$ & $1.086$ & $0.67$  \\
\hline\hline  
\end{tabular}
\end{table}

\begin{table}[t]
\caption{Accumulated cell ($\cal C$) and leaf ($\cal L$) counts in $10^6$ elements of adaptive 3D computations and resulting 
memory and mesh compression rates.}\label{tab:3DMemory}\centering\small 
\renewcommand{\arraystretch}{1.1}
\begin{tabular}{r||cc|cc||cc|cc}\hline\hline
\multirow{2}{*}{$L$} & \multicolumn{4}{c||}{MR} & \multicolumn{4}{c}{MRLT}\\ \cline{2-9}
& $\sum_n {\cal C}^n$ & [\%] & $\sum_n {\cal L}^n$ & [\%] & $\sum_n {\cal C}^n$ & [\%] & $\sum_n {\cal L}^n$ & [\%] \\ \hline
$5$ & $0.41$ & $78.0$ & $0.26$ & $48.8$ & $0.40$ & $77.1$ & $0.26$ & $49.0$ \\
$6$ & $3.39$ & $40.4$ & $2.03$ & $24.2$ & $3.32$ & $39.6$ & $2.00$ & $23.8$ \\
$7$ & $25.3$ & $18.9$ & $16.3$ & $12.1$ & $24.8$ & $18.5$ & $15.9$ & $11.9$ \\
$8$ &  $220$ & $10.2$ &  $147$ & $6.86$ &  $217$ & $10.1$ &  $146$ & $6.78$ \\
\hline\hline  
\multirow{2}{*}{$L$} & \multicolumn{4}{c||}{AMR} & \multicolumn{4}{c}{AMRLT}\\ \cline{2-9}
& $\sum_n {\cal C}^n$ & [\%] & $\sum_n {\cal L}^n$ & [\%] & $\sum_n {\cal C}^n$ & [\%] & $\sum_n {\cal L}^n$ & [\%] \\ \hline
$5$ & $0.34$ & $65.3$ & $0.26$ & $49.4$ & $0.34$ & $65.5$ & $0.27$ & $50.9$ \\
$6$ & $2.62$ & $31.3$ & $1.97$ & $23.5$ & $2.60$ & $30.9$ & $1.89$ & $22.6$ \\
$7$ & $22.1$ & $16.5$ & $16.6$ & $12.4$ & $21.9$ & $16.3$ & $16.1$ & $12.0$ \\
$8$ &  $220$ & $10.2$ &  $165$ & $7.68$ &  $219$ & $10.2$ &  $165$ & $7.68$ \\
\hline\hline  
\end{tabular}
\end{table}

The accuracy of the different schemes is analyzed in Table~\ref{tab:3Accuracy} by considering the $L_1$ error of density with respect to the reference solution, down-sampled to the corresponding level. The behavior of the $L_1$ errors is very similar
to the previous 2D study. For all cases, the error decreases as the number of levels $L$ is increased. Since again slightly 
different FV schemes are employed (cf. Section~\ref{sec:nummethods}), the MR/MRLT computations and their baseline FV scheme exhibit 
smaller absolute errors than the AMR/AMRLT methods and their respective FV scheme. For this configuration,
all methods show a convergence order being slightly above one at the finest level. We also observe that the adaptive 
computations with all four methods preserve the convergence order of the underlying FV scheme particularly well and 
even better than in the 2D investigation, which can be seen in Table~\ref{tab:3Accuracy} reflected in the very small error perturbation rates. 

Table~\ref{tab:3DMemory} gives the total cell and leaf counts over the time steps and the corresponding memory and mesh compression rates.
To keep the computational effort for the 3D case reasonable, the total numbers of used cells and leafs at the finest level $L=8$ is comparable to the 
2D values of Table~\ref{tab:2DMemory} at the respective finest level $L=10$. 
Comparing the cell and leaf counts in Table~\ref{tab:3DMemory} with 
Table~\ref{tab:2DMemory}, one finds that cell and leaf count increase roughly by a factor of 8 to 10 for each additional level in the 3D case, while factors $4$ to $5$ are 
found in the adaptive 2D simulations.  
This is consistent with an isotropic refinement by a factor of $2$ in the third dimension. Since the number of integrated cells in the FV simulations rises by the same factor, memory and mesh compression rates decrease with comparable factors in both the 
2D and the 3D configuration. 
However, in the 3D case the convergence is enhanced, with generally smaller compression ratios than in 2D at all levels except the coarsest one. 
Again, it can be noted that for MR/MRLT the compression rates decrease slightly faster than for the 
AMR schemes but the difference is generally smaller than in the 2D case. 

\begin{table}[t]
\caption{Computing times and CPU time compression rates for 3D computations.}\label{tab:3DCPU}\centering\small
\renewcommand{\arraystretch}{1.1}
\begin{tabular}{r||c|cc|cc||c|cc|cc}
\hline\hline
\multirow{2}{*}{$L$} & FV & \multicolumn{2}{c|}{MR} & \multicolumn{2}{c||}{MRLT}      & FV & \multicolumn{2}{c|}{AMR}  & 
\multicolumn{2}{c}{AMRLT}\\ \cline{3-6}\cline{8-11}
&[s] &[s] & [\%]&[s] & [\%]         &[s] &[s] & [\%]  & [s] & [\%]\\
 \cline{2-7}
\hline
$5$ & $7.495$ & $8.456$ &  $113$ & $8.414$ &  $112$ & $1.143$ & $1.250$ & $109$  & $1.038$ & $90.8$ \\
$6$ & $115.0$ & $67.29$ & $58.5$ & $65.34$ & $56.8$ & $16.79$ & $10.66$ & $63.5$ & $8.145$ & $48.5$ \\
$7$ & $1832$  & $509.0$ & $27.8$ & $492.0$ & $26.9$ & $260.9$ & $87.73$ & $33.6$ & $67.43$ & $25.8$ \\
$8$ & $28842$ & $4571$  & $15.8$ & $4477$  & $15.5$ & $4010$  & $863.5$ & $21.5$ & $676.7$ & $16.9$ \\
\hline\hline  
\end{tabular}
\end{table}

Denoting by $d$ the spatial dimension, the number of integrated cells in the FV simulations rises exactly by the factor $2^{d+1}$ at the next 
finer level. 
The absolute CPU times of Table~\ref{tab:3DCPU} exhibit the expected factor of $\sim 16$ increase for both FV schemes at all levels. 
In FV mode, the 3D MR/MRLT code is roughly a factor of $7.2$ slower than the 3D AMR/AMRLT code, which is an improvement compared to the 2D case. 
As in 2D, MR and MRLT methods show very little difference in run time, where the inefficiency of the AMR method \correction{without local time step} refinement, 
already observed in 2D, is actually clearly reduced in the 3D case. This behavior can be explained by the fact that the benefit of time step refinement can be 
expected to remain constant while the overall workload increases by a factor of $2$ with each additional spatial dimension.  
\correction{The overhead after Eq.~(\ref{eq:overhead}) for the AMR and AMRLT computations at $L=8$ is $180\%$ and $120\%$, respectively. For the MR and MRLT 
runs at this level, the overhead rate is $131\%$ and $129\%$ when being computed based on the performance of the MR code in FV unigrid code, and $1558\%$ and $1542\%$ when 
the AMR code in unigrid mode would be used as a somewhat questionable FV reference. See Section~\ref{sec:2dcomp} for a discussion.
Comparing the results with the 2D case of the previous section, one finds that the relative mesh adaptation overhead of the AMR/AMRLT computations rises considerably less 
than for MR/MRLT.} A consequence of this behavior is that the overhead in terms of the leaf updates is now consistently smaller for the AMRLT case than for MR but particularly 
also for MRLT. This result is different than in Section~\ref{sec:2dcomp}. However, for $L\ge 6$ the adaption overhead of the AMR method \correction{without local time stepping} 
is still larger than the overhead of MR or MRLT.

\begin{table}[t]
\caption{Breakdown of CPU time in \% spent in main task groups for 3D computations at $L=8$.}
\label{tab:Breakdown3D}\centering\small
\renewcommand{\arraystretch}{1.1}
\begin{tabular}{l|rrr||l|rrr}\hline\hline
\multirow{2}{*}{Task group} & \multicolumn{3}{c||}{MR} & \multirow{2}{*}{Task group} & \multicolumn{3}{c}{AMR} \\ \cline{2-4}\cline{6-8}
& \multicolumn{1}{c}{FV} & \multicolumn{1}{c}{NLT} & \multicolumn{1}{c||}{LT} & & \multicolumn{1}{c}{FV} & 
\multicolumn{1}{c}{NLT} & \multicolumn{1}{c}{LT} \\ \hline 
Numerics           & 34.44 & 23.43 & 23.23 & Numerics        & 92.12 & 69.51 & 74.04 \\ \cline{5-8}
Temp. data     & 38.99 & 16.21 & 16.03 & AMR data org.   &    -- &  7.84 &  5.52 \\ 
Boundary cond.     &  8.15 &  4.18 &  4.17 & Clustering      &    -- &  2.54 &  2.25 \\ \cline{1-4}
MR tree org.   &    -- & 34.48 & 34.56 & Flagging        &    -- &  2.64 &  1.61 \\
Level transfer     &    -- &  1.42 &  1.44 & Lever transfer  &    -- &  9.10 &  7.86 \\ \hline
Memory (c-lib)     & 12.61 &  5.77 &  5.69 & Memory (c-lib)  &  5.51 &  3.77 &  3.92 \\             
Unassigned         &  5.81 & 14.51 & 14.88 & Unassigned      &  2.37 &  4.60 &  4.80 \\ \hline \hline
\end{tabular}
\end{table}

An explanation for the enhanced performance of the AMRLT method might be found in Table~\ref{tab:Breakdown3D} that shows again a breakdown of the finest
resolved adaptive computation as obtained by the {\tt perf} tool. The breakdown in major task groups is the same as in Table~\ref{tab:breakdown}. 
Although the unassigned portion has substantially risen in the MR/MRLT cases, the time spent in the core FV update routine has shrunk slightly
while for AMR and AMRLT this portion has risen by about 5 and $2.5\,\%$, respectively. Further, the portion of the clustering algorithm has 
been reduced considerably in the 3D case. An explanation for this behavior might be that the flow phenomenon to be refined is more localized in the 3D simulations
and meshes are coarser. 
%

\subsection{Software design aspects}\label{sec:software}
The observed differences in performance warrant a deeper discussion of the internal
data structures and algorithmic solutions adopted in the computer codes which were employed for this study. 

\correction{The Carmen MR/MRLT code follows an object-oriented design whose base class is a {\it Cell} 
that stores the vector of state for a single FV cell.
When the code is operating in unigrid mode, a consecutive array of {\it Cell} elements is allocated; in  
MR/MRLT mode each {\it Cell} is a member of an object of type {\it Node}. In order to incorporate each {\it Node} object into the graded tree structure, {\it Node}
has a pointer to its parent and possibly an array with $2^d$ pointers to the children nodes. 
Depending on the dynamic mesh evolution, {\it Node} objects are dynamically 
allocated or deleted individually, using the standard \correction{\Cpp} commands {\tt new} and {\tt delete}. The children array remains unallocated for leaf nodes. 
When the numerical update is performed in Carmen, the nodes of the graded tree are recursively visited starting from a single root node. If the node is a leaf,
numerical fluxes are evaluated by querying the neighbor nodes for their cell-wise vector of state. All numerical fluxes along the cell boundary are then computed 
and their sum stored in {\it Cell}. This cell-oriented approach minimizes the memory footprint but requires two evaluations of the computationally expensive numerical 
fluxes per facet, even when the code is run in FV unigrid mode.
Since the graded tree is genuinely unstructured and the data are not stored 
consecutively in memory, the necessary tree traversal is computationally very costly as can be seen in Tables~\ref{tab:breakdown} and \ref{tab:Breakdown3D} 
(MR tree organization). 
Since only recursive tree traversal starting at the root node is currently supported, the execution of a higher level local time step requires a loop over 
all {\it Node} elements of the entire tree data structure, which explains the rather similar computing times of MR and MRLT.} 

\correction{The performance for accessing data through the tree in Carmen could be enhanced by implementing a recursive storage pattern that ensures data
locality in memory, \textit{e.g.}, by employing a generalized space filling curve algorithm to define a recursive sequential ordering of the multi-dimensional 
cells of the Cartesian mesh, cf. \cite{Weinzierl-Mehl-11}. The index information for the space filling curve as well as the local
neighbor and parent/child information of the tree nodes could be encoded very effectively in bit patterns, cf. \cite{Burstedde-etal-11}. Double flux computations could be relatively
easily avoided by adding available numerical fluxes at once to the vector of state of both neighboring cells or by introducing a more complex tree traversal operator.}

\correction{In the object-oriented design of block-based AMR software it is common to introduce a {\it Box} class defining a rectangular area in integer 
index space. A list of {\it Box} objects specifies the grid topology of every refinement level and a {\it Patch} class adds consecutive data storage to each {\it Box}. 
Depending on the stencil width of the numerical method, layers of ghost cells are created in addition. In AMROC, the allocation of 
{\it Patch} objects and their internal data arrays is also done with the standard \correction{\Cpp} {\tt new} and {\tt delete} commands but note 
that thanks to clustering typically several hundreds of cells together in a single patch, the number of memory management requests is much smaller than in Carmen. 
After the ghost cells of the {\it Patch} objects of one level have been properly set, the FV update function can be called in a simple loop 
over all patches. In unigrid mode, this corresponds to a straightforward block-based FV code on a rectangular grid and the computational performance hence is quasi-optimal.  
Each numerical flux is computed only once per patch; double flux evaluations occur only along the boundary facets that a patch shares with a neighbor of the same level. 
In explicit FV methods the computing costs are dominated by the numerical flux evaluations, which in combination with the faster MUSCL-Hancock reconstruction approach (cf. Section~\ref{sec:amr}) 
explains well why the absolute costs for the portion {\em Numerics} in Tables~\ref{tab:breakdown} and \ref{tab:Breakdown3D} 
are roughly 2 to 2.5 times smaller in the AMROC code.} 



%
%
%
%

\section{Conclusions}\label{conclusions}
In the present paper, fully adaptive computations of the 2D and 3D compressible Euler equations are presented
and two approaches for introducing adaptivity, \textit{i.e.}, MR and adaptive mesh refinement, are compared and benchmarked in detail.
For the same accuracy, the efficiency 
in terms of CPU time and memory
compression of these two adaptive methods has been assessed.
For time integration either global time stepping or scale-dependent local
time stepping techniques of second order Runge-Kutta type are used.
The main differences between MR and AMR is in the way the adaptive meshes
are stored and how the error estimators are defined.
The MR method uses a graded tree data structure and  thresholding of the
wavelet coefficients, corresponding to the details between two consecutive
levels, to define the adaptive mesh.
The AMR method uses a series of regular data blocks on the different levels
and relies on scaled gradient criteria based
on pressure and density to trigger the mesh refinement and coarsening.
Although in both approaches the thresholds are chosen to impose approximately the
same accuracy relative to a uniform grid computation, MR methods benefit from a rigorous and potentially more accurate
regularity analysis, while for
AMR methods rigorous error estimators are not available. Therefore,
threshold values of AMR have to be tuned
for a given problem, whereas in MR, in principle, the threshold is
independent of the problem.


\correction{The computational results show that the MR method generally presents larger compression rates
and has a prinicipal potential for obtaining larger gains in CPU time than the AMR method. The 
improved compression rates observed for the MR method are the consequences of a mathematically more sound adaptation criterion
as well as applying a patch-based refinement in the AMR method, which leads to a larger number of total cells to
avoid data fragmentation.}

When absolute CPU times are the primary concern, this investigation has underscored 
the importance of utilizing hierarchical data structures that preserve some memory coherence on computing data 
and use auxiliary data to avoid repeated generation of topological and numerical information. While the fulfillment of these requirements 
is rather naturally embedded into the block-based AMR approach, it comes at the cost of a more complicated base implementation.
On the other hand, the cell-based quad- or octree approach can be implemented prototypically with comparably little code; yet, the realization of 
a high-performance library for adaptive tree data is at least as involved, cf. \cite{Khokhlov-98,Weinzierl-Mehl-11,Burstedde-etal-11}, as devising
a high-performance AMR library, cf. \cite{Rendleman-etal-00,Hornung-etal-06,DeiterdingESAIM:2011}. 
\correction{In relative terms, the benefit of the MR approach has been clearly demonstrated. However, the herein employed prototype
code falls far short compared to the AMR approach in absolute compute time and -- with local time stepping -- even in 
relative adaption overhead for the 3D problem considered.}

\correction{Our next step will be to implement MR analysis as a mesh refinement criterion 
for the AMR/AMRLT algorithms, which will combine the natural high computational efficiency of the block-based approach with a
mathematically rooted error indicator, requiring less user adjustment. This paper has underscored that for typical, shock-dominated gas
dynamical problems solved with explicit FV schemes the choice of data structures and adaptation algorithms has a higher impact
on computational performance than the mathematical rigor of the mesh refinement criterion.}

\section*{Acknowledgments}
MOD and SG thankfully acknowledge financial support from
 Ecole Centrale de Marseille (ECM),
Funda\c{c}\~ao de
Amparo a Pesquisa do Estado de S\~ ao Paulo - FAPESP, and  the
Brazilian Research Council - CNPq, Brazil.
KS thanks the ANR project {}``SiCoMHD\char`\"{} for financial support. 
The authors also thank the CEMRACS 2012 summer program, 
where part of this study was done.
We are grateful to D. Foug\`ere and V. E. Mencone for their
computational assistance.


\begin{thebibliography}{10}

\bibitem{AbgrallHarten:1998}
R.~Abgrall and A.~Harten.
\newblock Multiresolution representation in unstructured meshes.
\newblock {\em SIAM J. Numer. Anal.}, 35(6):2128--2146, 1998.

\bibitem{BR01}
{R}. {B}ecker and {R}. {R}annacher.
\newblock An optimal control approach to a-posteriori error estimation.
\newblock In {A}. {I}serles, {R}. {B}ecker, {R}. {R}annacher, and {P}.~{G}.
  {C}iarlet, editors, {\em {A}cta {N}umerica}, volume~10, pages 1--102.
  {C}ambridge {U}niversity {P}ress, 2001.

\bibitem{Bell-Berger-Saltzman-94}
J.~Bell, M.~Berger, J.~Saltzman, and M.~Welcome.
\newblock Three-dimensional adaptive mesh refinement for hyperbolic
  conservation laws.
\newblock {\em {SIAM} J. Sci. Comput.}, 15(1):127--138, 1994.

\bibitem{Berger-Collela-88}
M.~Berger and P.~Colella.
\newblock Local adaptive mesh refinement for shock hydrodynamics.
\newblock {\em J. Comput. Phys.}, 82:64--84, 1988.

\bibitem{Berger-LeVeque-98}
M.~Berger and R.~J. LeVeque.
\newblock Adaptive mesh refinement using wave-propagation algorithms for
  hyperbolic systems.
\newblock {\em SIAM J. Numer. Anal.}, 35(6):2298--2316, 1998.

\bibitem{Berger-Oliger-84}
M.~Berger and J.~Oliger.
\newblock Adaptive mesh refinement for hyperbolic partial differential
  equations.
\newblock {\em J. Comput. Phys.}, 53:484--512, 1984.

\bibitem{Bihari:1996}
B.~L. Bihari.
\newblock Multiresolution schemes for conservation laws with viscosity.
\newblock {\em J. Comput. Phys.}, 123:207--225, 1997.

\bibitem{BihariHarten:1997}
B.~L. Bihari and A.~Harten.
\newblock Multiresolution schemes for the numerical solution of {2-D}
  conservation laws {I}.
\newblock {\em SIAM J. Sci. Comput.}, 18(2):315--354, 1996.

\bibitem{Brand-77}
A.~Brandt.
\newblock Multi-level adaptive solutions to boundary-value problems.
\newblock {\em Mathematics of Computations}, 31(183):333--390, April 1977.

\bibitem{Burstedde-etal-11}
C.~Burstedde, L.~C. Wilcox, and O.~Ghattas.
\newblock {\texttt{p4est}}: Scalable algorithms for parallel adaptive mesh
  refinement on forests of octrees.
\newblock {\em SIAM Journal on Scientific Computing}, 33(3):1103--1133, 2011.

\bibitem{ChiavassaDonat:2001}
G.~Chiavassa and R.~Donat.
\newblock Point value multi-scale algorithms for {2D} compressible flow.
\newblock {\em {SIAM} J. Sci. Comput.}, 23(3):805--823, 2001.

\bibitem{Cohen:2000}
A.~Cohen.
\newblock Wavelet methods in numerical analysis.
\newblock In P.~G. Ciarlet and J.~L. Lions, editors, {\em Handbook of
  {N}umerical {A}nalysis}, volume {VII}. Elsevier, Amsterdam, 2000.

\bibitem{CDKP:2000}
A.~Cohen, N.~Dyn, S.~M. Kaber, and M.~Postel.
\newblock Multiresolution finite volume schemes on triangles.
\newblock {\em J. Comput. Phys.}, 161:264--286, 2000.

\bibitem{CohenKaberMullerPostel:2003}
A.~Cohen, S.~M. Kaber, S.~M\"uller, and M.~Postel.
\newblock Fully adaptive multiresolution finite volume schemes for conservation
  laws.
\newblock {\em Math. Comp.}, 72:183--225, 2003.

\bibitem{DGM:2001}
W.~Dahmen, B.~Gottschlich-{M\"uller}, and S.~{M\"uller}.
\newblock Multiresolution schemes for conservation laws.
\newblock {\em Numer. Math.}, 88(3):399--443, 2001.

\bibitem{AMROC}
R.~Deiterding.
\newblock {AMROC - B}lockstructured {A}daptive {M}esh {R}efinement in
  {O}bject-oriented {C}++.
\newblock http://amroc.sourceforge.net.

\bibitem{Deiterding-PhDThesis}
R.~Deiterding.
\newblock {\em Parallel adaptive simulation of multi-dimensional detonation
  structures}.
\newblock PhD thesis, Brandenburgische Technische Universit\"at Cottbus, Sep
  2003.

\bibitem{Deiterding-03}
R.~Deiterding.
\newblock Construction and application of an {AMR} algorithm for distributed
  memory computers.
\newblock In T.~Plewa, T.~Linde, and V.~G. Weirs, editors, {\em Adaptive Mesh
  Refinement - Theory and Applications}, volume~41 of {\em Lecture Notes in
  Computational Science and Engineering}, pages 361--372. Springer, 2005.

\bibitem{Deiterding-08c}
R.~Deiterding.
\newblock A parallel adaptive method for simulating shock-induced combustion
  with detailed chemical kinetics in complex domains.
\newblock {\em Comp. Struct.}, 87:769--783, 2009.

\bibitem{DeiterdingESAIM:2011}
R.~Deiterding.
\newblock Block-structured adaptive mesh refinement - theory, implementation
  and application.
\newblock {\em {ESAIM} Proceedings}, 34:97--150, 2011.

\bibitem{DeiterdingDominguesGomesRousselSchneiderESIAM:2009}
R.~Deiterding, M.~O. Domingues, S.~M. Gomes, O.~Roussel, and K.~Schneider.
\newblock Adaptive multiresolution or adaptive mesh refinement? {A} case study
  for 2d {Euler} equations.
\newblock {\em {ESAIM} Proceedings}, 16:181--194, 2009.

\bibitem{VTF}
R.~Deiterding, R.~Radovitzki, S.~Mauch, F.~Cirak, D.~J. Hill, C.~Pantano, J.~C.
  Cummings, and D.~I. Meiron.
\newblock Virtual {T}est {F}acility: A virtual shock physics facility for
  simulating the dynamic response of materials.
\newblock http://www.vtf.website.

\bibitem{DominguesGomesDiaz:2003}
M.~O. Domingues, S.~M. Gomes, and L.~M.~A Diaz.
\newblock Adaptive wavelet representation and differenciation on
  block-structured grids.
\newblock {\em Appl. Num. Math.}, 47:421--437, 2003.

\bibitem{DGRS:2008}
M.~O. Domingues, S.~M. Gomes, O.~Roussel, and K.~Schneider.
\newblock An adaptive multiresolution scheme with local time stepping for
  evolutionary {PDEs}.
\newblock {\em {J. Comp. Phys.}}, 227:3758--3780, 2008.

\bibitem{DGRS:2009}
M.~O. Domingues, S.~M. Gomes, O.~Roussel, and K.~Schneider.
\newblock Space-time adaptive multiresolution methods for hyperbolic
  conservation laws: {A}pplications to compressible {Euler} equations.
\newblock {\em Appl. Num. Math.}, 59:2303--2321, 2009.

\bibitem{DominguesGomesRousselSchneiderESAIM:2011}
M.~O. Domingues, S.~M. Gomes, O.~Roussel, and K.~Schneider.
\newblock Adaptive multiresolution methods.
\newblock {\em {ESAIM} Proceedings}, 34:1--96, 2011.

\bibitem{DRS:2009}
M.~O. Domingues, O.~Roussel, and K.~Schneider.
\newblock An adaptive multiresolution method for parabolic pdes with time-step
  control.
\newblock {\em {Int. J. Numer. Meth. Engng.}}, 78:652--670, 2009.

\bibitem{Friedel-Grauer-Marliani-97}
H.~Friedel, R.~Grauer, and C.~Marliani.
\newblock Adaptive mesh refinement for singular current sheets in
  incompressible magnetohydrodynamics flows.
\newblock {\em J. Comput. Phys.}, 134(1):190--198, 1997.

\bibitem{Gottschlich:1999}
B.~Gottschlich-M\"uller and S.~M\"uller.
\newblock {Adaptive finite volume schemes for conservation laws based on local
  multiresolution techniques}.
\newblock In R.~Jeltsch and M.~Frey, editors, {\em Hyperbolic Problems: Theory,
  Numerics, Applications}, volume {129}. ISNM, Inter. Ser. Numer. Math., 1999.

\bibitem{Harten:1995}
A.~Harten.
\newblock Multiresolution algorithms for the numerical solution of hyperbolic
  conservation laws.
\newblock {\em Comm. Pure Appl. Math.}, 48:1305--1342, 1995.

\bibitem{Harten:1996}
A.~Harten.
\newblock Multiresolution representation of data: a general framework.
\newblock {\em {SIAM} {J. Numer. Anal.}}, 33(3):385--394, 1996.

\bibitem{Holmstrom:1997}
M.~Holmstr{\"o}m.
\newblock {\em {Wavelet Based Methods for Time Dependent {PDEs}}}.
\newblock PhD thesis, Uppsala University, 1997.

\bibitem{Holmstrom:1999}
M.~Holmstr\"om.
\newblock Solving hyperbolic {PDEs} using interpolating wavelets.
\newblock {\em SIAM J. Sci. Comput.}, 21(2):405--420, 1999.

\bibitem{Hornung-etal-06}
R.~D. Hornung, A.~M. Wissink, and S.~H. Kohn.
\newblock Managing complex data and geometry in parallel structured {AMR}
  applications.
\newblock {\em Engineering with Computers}, 22:181--195, 2006.

\bibitem{Kaibara:2000}
M.~Kaibara and S.~M. Gomes.
\newblock {A fully adaptive multiresolution scheme for shock computations}.
\newblock In E.F. Toro, editor, {\em Godunov Methods: Theory and Applications},
  volume~{}. Klumer Academic/Plenum Publishers, 2001.

\bibitem{Khokhlov-98}
A.~M. Khokhlov.
\newblock Fully threaded tree algorithms for adaptive refinement fluid dynamics
  simulations.
\newblock {\em J. Comput. Phys.}, 143:519--543, 1998.

\bibitem{Kohn-Baden-95}
S.~R. Kohn and S.~B. Baden.
\newblock A parallel software infrastructure for structured adaptive mesh
  methods.
\newblock In {\em Proc. of the Conf. on Supercomputing '95}, December 1995.

\bibitem{Kolomenskiy:2015}
D.~Kolomenskiy, J.-C. Nave, and K.~Schneider.
\newblock Adaptive gradient-augmented level set method with multiresolution
  error estimation.
\newblock {\em J. Sci. Comput.}, 2015.
\newblock doi:10.1007/s10915-015-0014-7, in press.

\bibitem{LaxLiu:1998}
P.~D. Lax and X.~D. Liu.
\newblock Solution of two-dimensional {Riemann problems} of gas dynamics by
  positive schemes.
\newblock {\em {SIAM J. Sci. Comput.}}, 19(2):319--340, 1998.

\bibitem{Liou:1996}
M.-S. Liou.
\newblock A sequel to {AUSM}: {AUSM+}.
\newblock {\em J. Comput. Phys.}, 129:364--382, 1996.

\bibitem{MacNeice-etal-00}
P.~MacNeice, K.~M. Olson, C.~Mobarry, R.~deFainchtein, and C.~Packer.
\newblock {PARAMESH: A} parallel adaptive mesh refinement community toolkit.
\newblock {\em Computer Physics Communications}, 126:330--354, 2000.

\bibitem{Muller:2003}
S.~M\"uller.
\newblock {\em Adaptive multiscale schemes for conservation laws}, volume~27 of
  {\em Lectures Notes in Computational Science and Engineering}.
\newblock Springer, Heidelberg, 2003.

\bibitem{MullerStiriba:2007}
S.~M\"uller and Y.~Stiriba.
\newblock Fully adaptive multiscale schemes for conservation laws employing
  locally varying time stepping.
\newblock {\em J. Sci. Comput.}, 30(3):493--531, 2007.

\bibitem{Pantano-Deiterding-05}
C.~Pantano, R.~Deiterding, D.~J. Hill, and D.~I. Pullin.
\newblock A low-numerical dissipation patch-based adaptive mesh refinement
  method for large-eddy simulation of compressible flows.
\newblock {\em J. Comput. Phys.}, 221:63--87, 2007.

\bibitem{Rendleman-etal-00}
C.~A. Rendleman, V.~E. Beckner, M.~Lijewski, W.~Crutchfield, and J~B. Bell.
\newblock Parallelization of structured, hierarchical adaptive mesh refinement
  algorithms.
\newblock {\em Computing and Visualization in Science}, 3, 2000.

\bibitem{RossinelliHejazialhosseiniSpampinatoKoumoutsakos:2011}
D.~Rossinelli, B.~Hejazialhosseini, D.~Spampinato, and P.~Koumoutsakos.
\newblock Multicore/multi-gpu accelerated simulations of multiphase
  compressible flows using wavelet adapted grids.
\newblock {\em SIAM J. Sci. Comput.}, 33:512--540, 2011.

\bibitem{Rossinellietal:2015}
D.~Rossinelli, B.~Hejazialhosseini, W.~van Rees, M.~Gazzola, M.~Bergdorf, and
  P.~Koumoutsakos.
\newblock Mrag-i2d: Multiresolution adapted grids for remeshed vortex methods
  on multicore architectures.
\newblock {\em J. Comput. Phys.}, 288:1--18, 2015.

\bibitem{RousselSchneider:2005}
O.~Roussel and K.~Schneider.
\newblock An adaptive multiresolution method for combustion problems:
  application to flame ball - vortex interaction.
\newblock {\em Comp. Fluids}, 34(7):817--831, 2005.

\bibitem{RSTB03}
O.~Roussel, K.~Schneider, A.~Tsigulin, and H.~Bockhorn.
\newblock A conservative fully adaptive multiresolution algorithm for parabolic
  {PDEs}.
\newblock {\em J. Comput. Phys.}, 188(2):493--523, 2003.

\bibitem{SchneiderVasilyev:2010}
K.~Schneider and O.~V. Vasilyev.
\newblock Wavelet methods in computational fluid dynamics.
\newblock {\em Annu. Rev. Fluid. Mech.}, 42:473--503, 2010.

\bibitem{SchulzRinneetal:1993}
C.~W. Schulz-Rinne, J.~P. Collis, and H.~M. Glaz.
\newblock Numerical solution of the {Riemann} problem for two-dimensional gas
  dynamics.
\newblock {\em {SIAM} J. Sci. Comput.}, 14:1394--1414, 1993.

\bibitem{Toro}
E.~F. Toro.
\newblock {\em Riemann solvers and numerical methods for fluid dynamics}.
\newblock Springer, 1997.

\bibitem{Wada-Liou-97}
Y.~Wada and M.~S. Liou.
\newblock An accurate and robust flux splitting scheme for shock and contact
  discontinuities.
\newblock {\em {SIAM J. Sci. Comput.}}, 18(3):633--657, 1997.

\bibitem{Weinzierl-Mehl-11}
T.~Weinzierl and M.~Mehl.
\newblock Peano -- {A} traversal and storage scheme for octree-like adaptive
  cartesian multiscale grids.
\newblock {\em SIAM Journal on Scientific Computing}, 33(5):2732--2760, 2011.

\bibitem{ZhangZheng:1990}
T.~Zhang and Y.~Zheng.
\newblock Conjecture on the structure of solutions the {Riemann} problem for
  two-dimensional gas dynamics systems.
\newblock {\em {SIAM} J. Math. Anal.}, 21:593--630, 1990.

\end{thebibliography}

\end{document}